\documentclass[10pt,a4paper]{article}
\usepackage[draft]{optional}
\usepackage[british,american]{babel}
\usepackage{latexsym}
\usepackage{amsxtra}
\usepackage{mathrsfs}
\usepackage{amsfonts, amsmath, wasysym}
\usepackage{amscd}
\usepackage{slashed}
\usepackage{enumerate}
\usepackage{amssymb,amsthm,
paralist
}

\usepackage{
latexsym,
nicefrac,
}

\usepackage[usenames]{color}

\usepackage{url}

\definecolor{darkgreen}{rgb}{0,0.5,0}
\definecolor{darkred}{rgb}{0.7,0,0}
\usepackage[colorlinks,
citecolor=darkgreen, linkcolor=darkred
]{hyperref}
\usepackage{esint}
\usepackage{bibgerm}
\usepackage[normalem]{ulem}

\theoremstyle{plain}
\newtheorem{lemma}{Lemma}[section]
\newtheorem{thm}[lemma]{Theorem}
\newtheorem{prop}[lemma]{Proposition}
\newtheorem{cor}[lemma]{Corollary}
\theoremstyle{definition}
\newtheorem{defi}[lemma]{Definition}

\newtheorem{rmk}[lemma]{Remark}

\numberwithin{equation}{section}

\newcommand{\beq}{\begin{equation}}
\newcommand{\eeq}{\end{equation}}
\newcommand{\beqa}{\begin{equation}\begin{aligned}}
\newcommand{\eeqa}{\end{aligned}\end{equation}}
\newcommand{\brmk}{\begin{rmk}}
\newcommand{\ermk}{\end{rmk}}
\newcommand{\partref}[1]{\hbox{(\csname @roman\endcsname{\ref{#1}})}}

\def\bee{\begin{eqnarray}}
\def\beee{\begin{eqnarray*}}
\def\eee{\end{eqnarray}}
\def\eeee{\end{eqnarray*}}
\def\nn{\nonumber}
\def\ba{\begin{array}}
\def\ea{\end{array}}

\newcommand{\Dir}{D \hskip -2mm \slash}
\newcommand{\Par}{\partial \hskip -1.8mm \slash}

\def\la{\langle}
\def\ra{\rangle}

\def\RR{\Bbb R}

\def \a {\alpha}
\def \b {\beta}

\def \g {\gamma}

\def \ra {\longrightarrow}

\def\la{\langle}
\def\ra{\rangle}
\def\cal{\mathcal}

\newcommand{\twomat}[4]{ \left( \begin{array}{cc}
#1 & #2 \\
#3 & #4
\end{array} \right)}

\newcommand{\cvec}[2] {\left( \begin{array}{c}
#1 \\
#2
\end{array}\right)}

\newcommand{\twopartdef}[4]
{
	\left\{
		\begin{array}{ll}
			#1 & \mbox{if  } #2 \bigskip \\
			#3 & \mbox{if  } #4
		\end{array}
	\right.
}

\newcommand{\pl}[2]{{\frac{\partial #1}{\partial #2}}}

\newcommand{\zb}{\overline{z}}
\newcommand{\db}{\overline{\partial}}

\newcommand{\R}{\mathbb{R}}

\newcommand{\Ca}{\mathcal{C}}
\newcommand{\M}{M}

\newcommand{\emb}{\hookrightarrow}
\newcommand{\cemb}{\subset \subset}
\newcommand{\In}{\subset}
\newcommand{\Om}{\Omega}
\newcommand{\om}{\omega}
\newcommand{\dl}{{\delta}}
\newcommand{\Dl}{{\Delta}}

\newcommand{\al}{{\alpha}}

\newcommand{\ed}{{\rm d}}

\newcommand{\dv}{{\rm d}^{\ast}}
\newcommand{\D}{{\nabla}}

\newcommand{\ti}[1]{{\tilde{#1}}}
\newcommand{\fr}[2]{\frac{#1}{#2}}
\newcommand{\sm}{{\setminus}}

\newcommand{\T}{{\rm T}}

\def\XXint#1#2#3{{\setbox0=\hbox{$#1{#2#3}{\int}$}
\vcenter{\hbox{$#2#3$}}\kern-.5\wd0}}

\title{\sc{Regularity at the free boundary for Dirac-harmonic maps from surfaces}
\\
}
\author{Ben Sharp\footnote{ \noindent Imperial College London, Huxley Building, Department of Mathematics, 180 Queen's Gate, SW7 2AZ, UK. \newline e-mail: ben.g.sharp@gmail.com}\,
and  Miaomiao Zhu\footnote{Max Planck Institute for Mathematics in the Sciences, Inselstrasse 22, 04103 Leipzig, Germany. \newline e-mail: Miaomiao.Zhu@mis.mpg.de }}

\date{\today}

\begin{document}

\maketitle

\begin{abstract}
We establish the regularity theory for certain critical elliptic systems with an anti-symmetric structure under inhomogeneous Neumann and Dirichlet boundary constraints.
As applications, we prove full regularity and smooth estimates at the free boundary for weakly Dirac-harmonic maps from spin Riemann surfaces.
Our methods also lead to the full interior $\epsilon$-regularity and smooth estimates for weakly Dirac-harmonic maps in all dimensions.
\end{abstract}

\section{Introduction}

Motivated by the supersymmetric nonlinear sigma model from quantum field theory \cite{De, J3}, a variational problem that couples a map between
two Riemannian manifolds with a spinor field along this map is introduced in \cite{CJLW1}. More precisely, let $(M,g_M)$ be a spin Riemannian manifold of dimension
$m \geq 2$,  $\Sigma M$ the spinor bundle over $M$ and $(N, g)$ a compact Riemannian manifold of dimension $d$. Let $\phi$ be a map from $M$ to $N$ and $\psi$ a section of the twisted bundle $\Sigma M \otimes\phi^{-1}TN$. Using the connection $\widetilde{\nabla}$ which is induced from the spin connection on
$\Sigma M$ and the pull back Levi-Civita connection on $\phi^{-1}TN$, one defines the Dirac operator $\Dir$ along the map $\phi$ by $\Dir \psi:=\g_\a\cdot\widetilde{\nabla}_{\g_\a}\psi$,
where $\g_\a$ is a local orthonormal frame on $M$ and $``\cdot"$ is the Clifford multiplication. Consider the functional
\bee \label{L} L(\phi,\psi) =\int_M \left ( |d\phi|^2+\langle \psi,\Dir \psi \rangle \right ) dV_{g_M}.
\eee
Critical points $(\phi,\psi)$ of \eqref{L} are called Dirac-harmonic maps from $M$ to $N$. Mathematically, Dirac-harmonic maps generalise the notions of harmonic
maps and harmonic spinors, both of which have been extensively studied from various geometric and analytic points of view.
Analytically, the Euler-Lagrange system of \eqref{L} couples a second order critical elliptic equation with a first order Dirac-type equation.
In dimension $m=2$, the functional \eqref{L} is conformally invariant in the domain as is the case for harmonic maps, which places the study of
Dirac-harmonic maps in the framework of two dimensional conformally invariant variational problems, where a lot of powerful geometric analysis techniques
have been developed, allowing for a deep investigation of the regularity theory for weak solutions.

Exploring the geometric and analytic aspects of this variational problem, Chen-Jost-Wang-Zhu \cite{CJWZ} introduced an
appropriate boundary value problem for Dirac-harmonic maps from a spin Riemann surface $M$ with non-empty boundary $\partial M$ into a
compact Riemannian manifold $N$, providing a mathematical interpretation of the supersymmetric nonlinear sigma model with boundaries \cite{ALZ1}
and the D-branes in superstring theory \cite{Po}. More precisely, let $\mathcal{S}$ be a closed $s$-dimensional submanifold of $N$, then
the map $\phi$ is supposed to satisfy the classical free boundary condition: $\phi(\partial M)\In \mathcal{S},$ while the spinor field $\psi \in \Gamma(\Sigma M\otimes \phi^{-1}TN)$ is required to satisfy a chirality type (local elliptic) boundary condition that is
compatible with $\phi$ as well as the supporting submanifold $\mathcal{S}$:
\bee\label{chir}
\mathbf{B}^{\pm}_{\phi}\psi|_{\partial M} = 0, \eee
generalising the usual chirality (local elliptic) boundary condition $\mathbf{B}^{\pm} \ti\psi|_{\partial M} = 0$ introduced by Gibbons-Hawking-Horowitz-Perry
\cite{GHHP} (we refer to \cite{HMR} for the mathematical setting) for usual spinors $\ti\psi\in \Gamma(\Sigma M)$.

Embed $N$ isometrically into some Euclidean space $\mathbb{R}^n$ and consider admissible fields for \eqref{L} in the following space
\bee   \mathcal{X}^{1,2}_{1,\fr43}(M,N;\mathcal{S}) := \left \{
\ba{ll}
 (\phi,\psi)  \in H^1(M,N) \times W^{1,\fr43}(\Sigma M \otimes \phi^{-1}TN): \nn\\ \phi(x) \in \mathcal{S} \  \text{and}  \
 \mathbf{B}_{\phi}^{\pm}\ \psi (x) = 0  \text{ for  a.e. } x \in  \partial M \nn
\ea
\right \}.
\eee
One can easily verify that the functional \eqref{L} extends to $\mathcal{X}^{1,2}_{1,\fr43}(M,N;\mathcal{S})$.
Critical points $(\phi,\psi)\in \mathcal{X}^{1,2}_{1,\fr43}(M,N; \mathcal{S})$ of \eqref{L} are called weakly Dirac-harmonic maps from $M$ to $N$
with free boundary $\phi(\partial M)$ on $\mathcal{S}$.

To study the regularity at the free boundary for weakly Dirac-harmonic maps from spin Riemann surfaces, by conformal invariance in dimension $m=2$, one can locate the problem in a small neighbourhood of a boundary point and consider the model case that the domain
is the upper half unit disc $$B^+_1:= \left \{(x_1,x_2) \in \mathbb{R}^2 | x_1^2+x_2^2 < 1, x_2 \geq 0  \right \}$$equipped with the Euclidean metric and the free boundary
portion is $I:= \left \{(x_1,0) \in \mathbb{R}^2 | -1 < x_1 < 1 \right \}.$ Denote by $\overrightarrow{n}$ the outward unit normal vector field on the boundary portion $I$.

As in \cite{CJWZ}, applying a reflection procedure for $\phi$ analogous to that used in the free boundary problem for harmonic maps \cite{GJ, Sc} and introducing
a reflection for $\psi$ that is compatible with the boundary condition \eqref{chir}, one can extend the two fields across $I$ to the whole ball $B_1$ such that the extended fields $(\phi,\psi)$
weakly solve a critical elliptic system with an anti-symmetric structure similar to that introduced by Rivi\`ere \cite{Ri} and Rivi\`ere-Struwe \cite{riviere_struwe} up to an additional frame transformation which can be analytically well controlled.
By adapting the regularity theory of Rivi\`ere-Struwe, they proved the following regularity results:
\begin{thm}[\cite{CJWZ}, Theorem 1.1-1.2] \label{oldthm} Let $M$ be a compact spin Riemann surface with boundary $\partial M$,
$N$ a compact Riemannian manifold with $\mathcal{S}\In N$ a closed totally geodesic submanifold of $N$.
Let $(\phi,\psi)$ be a weakly Dirac-harmonic map from $M$ to $N$ with free boundary $\phi(\partial M)$ on $\mathcal{S}$.
Then there exists some $\b \in (0,1)$ such that $ \phi \in C^{1,\b}(M,N)$ and $\psi \in C^{1,\b}(\Sigma M \otimes\phi^{-1}TN).$
\end{thm}
In the present paper, we prove the full regularity at the free boundary - without assuming that $\mathcal{S}\In N$ is totally geodesic -
by establishing the regularity theory for critical elliptic systems with an anti-symmetric structure under inhomogeneous Neumann and Dirichlet boundary constraints,
extending the classical boundary regularity theory for sub-critical elliptic PDE studied by Agmon-Douglis-Nirenberg \cite{agmon_doug_niren}. It  provides the boundary regularity results related to
the recent interior regularity theory for critical elliptic systems with an anti-symmetric structure developed by Rivi\`ere \cite{Ri},
Rivi\`ere-Struwe \cite{riviere_struwe}, Sharp-Topping \cite{Sh_To}, Sharp \cite{Sh, thesis} and Zhu \cite{Zhu}.

To state our main PDE result, let $ 1 \leq k \in \Bbb{N} $, $1 \leq p \in \Bbb{R}$ and $U\In \R^m$, $m\geq 2$ be some open set with $T\In \partial U$ a smooth boundary portion, define the space (see e.g. \cite{wehrheim})
$$W_{\partial}^{k,p}(T):=\{g\in L^1(T) : \text{$g = G|_{T}$ for some $G\in W^{k,p}(U)$}\}$$ with norm $\|g\|_{W_{\partial}^{k,p}(T)} = \inf_{G\in W^{k,p}(U), G|_T =g} \|G\|_{W^{k,p}(U)}.$
Then we have

\begin{thm}\label{thm: mainPDE}
Let $d \geq2$, $0\leq s\leq d$, $0<\Lambda <\infty$ and $1<p<2$. For any $A\in L^{\infty}\cap W^{1,2}(B^+_1, GL(d))$,
$\Omega\in L^2(B^+_1,so(d)\otimes \wedge^1 \R^m)$, $f\in L^p(B^+_1,\R^d)$, $(\frak{g},\frak{k})\in (W_{\partial}^{1,p}(I,\R^s),W_{\partial}^{2,p}(I,\R^{d-s})) $
and for any $u\in W^{1,2}(B^+_1,\R^d)$ weakly solving
\begin{eqnarray}
\dv(A\ed u) &=& \langle \Omega , A \ed u\rangle  + f \quad \text{in}\quad B_1^+,  \label{eqn: mainPDE}  \\
\pl{u^i}{\overrightarrow{n}} &=&\frak{g}^i \quad \text{on $I$}\quad 1\leq i \leq s,  \label{eqn: bcmain} \\
u^j &=& \frak{k}^{j} \quad \text{on $I$}\quad s+1 \leq j \leq d,     \label{eqn: bcmain1}
\end{eqnarray}
with
\begin{equation*}\label{eqn: Abound}
\Lambda^{-1}|\xi| \leq |A(x)\xi| \leq \Lambda |\xi| \,\,\,\,\,\text{for a.e. $x\in B_1^+$}, \,\,\,\text{for all $\xi\in \R^d$}
\end{equation*}
and satisfying the compatibility condition that $A(x)$ commutes with $R$ for a.e. $x\in B_1^+$ where
\begin{equation} \label{R}
R=(R^i_j):=\twomat{Id_s}{0}{0}{-Id_{d-s}}.
\end{equation}
Then $u\in W^{2,p}_{loc}(B_1^+\cup I)$ and there exist $\epsilon=\epsilon(d,\Lambda, p)>0$ and $C=C(d,\Lambda, p)>0$ such that whenever
$\|\Om\|_{L^2(B^+_1)}+
\|\D A\|_{L^2(B^+_1)}\leq \epsilon $ then
$$\|\D^2 u\|_{L^p(B^+_{\fr12})} + \|\D u\|_{L^{\fr{2p}{2-p}}(B^+_{\fr12})} \leq C(\| u\|_{L^1(B_1^+)} + \|f\|_{L^p(B_1^+)} +\|\frak{g}\|_{W_{\partial}^{1,p}(I)} +\|\frak{k}\|_{W_{\partial}^{2,p}(I)}).$$

\end{thm}
\begin{rmk}
Our proof of Theorem \ref{thm: mainPDE} does not follow from the `perturbed Coulomb frame' approach of Rivi\`ere \cite{Ri} (see also \cite{Sh_To}) but first requires the methods of Rivi\`ere-Struwe \cite{riviere_struwe} in obtaining an interior H\"older estimate (see \cite{Zhu} or Theorem \ref{thm: holder}). Then one uses the work of Sharp \cite{Sh} in order to bootstrap this into integrability of the gradient above $L^2$ (see Theorem \ref{thm: main2}). The advantage of this approach is that it works in all dimensions $m\geq 2$ in obtaining higher integrability on the interior. In two dimensions this essentially follows from the stability of the local Hardy space $h^1$ under multiplication by H\"older continuous functions, coupled with the embedding $h^1\emb H^{-1}$. To obtain sub-critical integrability in higher dimensions (i.e. in order to prove Theorem \ref{thm: main2} for $m>2$) one needs to extend this embedding and show $h^1\cap M^{1,m-2}\emb H^{-1}$ - see \cite{Sh} for details. The fundamental idea behind obtaining the regularity for such critical equations with anti-symmetric structures began with the work of Rivi\`ere \cite{Ri} and later Rivi\`ere-Struwe \cite{riviere_struwe}, which involved re-writing the equation with respect to a Coulomb-type frame, generalising the work of H\'elein \cite{helein_conservation}, thus uncovering Jacobian determinant structures which allows one to use Hardy space methods via the pioneering work of Coiffman, Lions, Mayer and Semmes \cite{clms}.
\end{rmk}

\begin{rmk}
Theorem \ref{thm: mainPDE} says that for the boundary conditions above, under the appropriate assumptions on $A$ and $\Om$, the solution $u$
behaves like that of the following classical boundary value problem (for a scalar-valued $u$)
$$-\Dl u = f  \quad \text{in $ B^+_1 $}  $$
$$\pl{u}{\overrightarrow{n}} = \frak{g} \quad \text{on $I$} \quad \quad\text{or} \quad\quad u = \frak{k} \quad \text{on $I$}$$
studied by Agmon-Douglis-Nirenberg \cite{agmon_doug_niren}, where
$f\in L^p(B^+_1)$, $\frak{g} \in W_{\partial}^{1,p}(I), \frak{k} \in W_{\partial}^{2,p}(I)$ and $p>1$.
\end{rmk}

%

\begin{rmk}
When the boundary conditions are homogeneous, i.e., $\frak{g}\equiv 0$ and $\frak{k}\equiv 0$, we have an extension of  Theorem \ref{thm: mainPDE} to the case of
higher dimensional domains $B_1^+\subset \R^m$ for $m\geq 2$  with $\Om$, $\D A$ and $\D u$ in the Morrey space $M^{2,m-2}$ and $f\in L^p$ with $\fr{m}{2}<p<m$
(see Section \ref{impPDE}, Theorem \ref{thm: mainPDEhom}).
\end{rmk}

Now, going back to the free boundary problem for Dirac-harmonic maps, it is possible to use the known interior regularity \cite{CJWZ, WX} to prove that we can localise our problem also in the
target - see Lemma \ref{nhdlemma}.  This allows us to use the so called {\it Fermi coordinates} about some point $q\in\mathcal{S}$ at which point we use elementary geometric
arguments to show that our map $\phi$ weakly solves an elliptic system of the form \eqref{eqn: mainPDE} with boundary conditions \eqref{eqn: bcmain} and \eqref{eqn: bcmain1}
but only for $p=1$. More precisely we have
$$ \text{$|\Omega(x)| + |\nabla A(x)| \leq C |\nabla \phi(x)|$,\, and \,$|f(x)| \leq C|\psi(x)|^2 |\D \phi(x)|$ for \textit{a.e.} $x \in  B^+_1$.}$$ Moreover the tangent components of the map $\phi$ satisfy an inhomogeneous Neumann condition:
$$\left(\pl{\phi^i}{\overrightarrow{n}}\right) = {\rm Re} (\mathcal{P}^i_{\mathcal{S}}(\overrightarrow{n}\cdot\psi^{\bot};\psi^{\top}) )\,\,\,\,\,\,\text{on $I$ for $1\leq i \leq s$}$$ and the normal components satisfy a
homogeneous Dirichlet condition
$$\text{$\phi^{i}|_I = 0$ for $ s+1 \leq i \leq d.$}$$ Here $\mathcal{P}_{\mathcal{S}}(\cdot; \cdot)$ is a spinorial extension of the shape operator
$P_{\mathcal{S}}(\cdot; \cdot)$ of $\mathcal{S}$ in $N$, see section \ref{recapfbvp}.
{\emph{A-priori} $\psi\in W^{1,\fr43}\emb L^4$ and $\phi \in W^{1,2}$ hence we have
$$\text{$f\in L^1(B^+_1)$ and ${\rm Re} (\mathcal{P}^i_{\mathcal{S}}(\overrightarrow{n}\cdot\psi^{\bot};\psi^{\top}) ) \in W_{\partial}^{1,1}(I).$}$$

In order to be able to apply the regularity result in Theorem \ref{thm: mainPDE}, we need to improve the regularity up to the boundary for the spinor field $\psi$
so that
\begin{equation}\label{www}
\text{$f \in L^p(B^+_1)$ and ${\rm Re} (\mathcal{P}^i_{\mathcal{S}}(\overrightarrow{n}\cdot\psi^{\bot};\psi^{\top}) )  \in W_{\partial}^{1,p}(I)$, for some $p>1$.}
\end{equation}
Fortunately, this can be achieved by investigating the Dirac equation for $\psi$ as well as observing that the chirality boundary condition \eqref{chir} on the spinor is morally a Riemann-Hilbert boundary condition
for the Cauchy-Riemann operator (see Theorems \ref{spinor_thm2} and \ref{spinor_thm}) and hence Theorem \ref{spinor_thm1} allows us to consider our map $\phi$ as above but this time \eqref{www} is satisfied for some $1<p<2$ and Theorem \ref{thm: mainPDE} yields, in particular, that $\D \phi \in L^r$ up to the boundary for some $r>2$.

We note that if we are only able to prove H\"older regularity of $\phi$ then the argument does not work - thus Theorem \ref{thm: mainPDE} is integral in obtaining the full regularity. The H\"older regularity is enough if one assumes the submanifold $\mathcal{S}$ is totally geodesic, since the tangent components of the map $\phi$ then satisfy a homogeneous Neumann condition and the usual reflection principle applies as in the case for Harmonic maps - see \cite{Sc}.
At this point we are considering a sub-critical and classical boundary value problem and a simple bootstrapping argument concludes:

\begin{thm}\label{thm: breg}
Let $N$ be a compact Riemannian manifold, $\mathcal{S}\In N$ a closed submanifold with $(\phi,\psi)$ a weakly Dirac-harmonic map from $B^+_1$ to $N$
 with free boundary $\phi(I)$ on $\mathcal{S}$. Then for any $k\in \mathbb{N}$ there exist
$\epsilon = \epsilon(N,\mathcal{S})>0$ and $C=C(k,N,\mathcal{S})>0$ such that if
$\| \D \phi \|_{L^2(B_1^+)} + \|\psi\|_{L^4(B_1^+)}
\leq \epsilon$, then
\bee \|\D^k \phi\|_{L^{\infty}(B_{\fr12}^+)} +\|\D^k \psi\|_{L^{\infty}(B_{\fr12}^+)} \leq C\| \D \phi \|_{L^2(B_1^+)}(1 + \| \psi \|_{L^4(B_1^+)} ). \nn
\eee
\end{thm}

By the conformal invariance of the problem and re-scaling in the domain coupled with the known interior regularity \cite{CJWZ, WX},
we trivially establish the full regularity at the free boundary:

\begin{thm}\label{newthm} Let $M$ be a compact spin Riemann surface with boundary $\partial M$,
$N$ a compact Riemannian manifold, and $\mathcal{S}\In N$ a closed submanifold.
Let $(\phi,\psi)$ be a weakly Dirac-harmonic map from $M$ to $N$ with free boundary $\phi(\partial M)$ on $\mathcal{S}$.
Then $(\phi,\psi)$ is smooth up to the boundary. \end{thm}

\begin{rmk}
We point out that we expect to be able to extend this result to a partial regularity result for weakly stationary Dirac-harmonic maps with a free boundary constraint in any dimension $m\geq 3$. This however will require an extra assumption on the domain  $M$ in order to make sense of the Chirality-type boundary condition. Moreover one would need a version of Theorem \ref{thm: mainPDE} for higher dimensions involving new techniques which are not necessary in this paper. The authors will deal with these issues in a subsequent work.
\end{rmk}

\begin{rmk}
Any weakly harmonic map or weakly minimal surface $u:B_1^+ \to N$ with a free boundary $u(I) \In\mathcal{S}$ is also a weakly Dirac-harmonic map when coupled with the zero spinor field.
Therefore we can apply Theorems \ref{thm: breg} and \ref{newthm} to obtain the regularity at the boundary in this case.
The free boundary value problem is simpler for harmonic maps and minimal surfaces since one has zero Neumann conditions for the tangential components as well as
zero Dirichlet conditions for the normal components - geometrically this is known as the transversality property at the free boundary (see e.g. \cite{Hil, GJ}):
$\pl{u}{\overrightarrow{n}} \bot T_u \mathcal{S},  a.e.  \text{ on} \  I$.
\end{rmk}

Using the PDE results from this paper we also obtain the following local smooth interior estimates for weakly Dirac-harmonic maps in all dimensions $m\geq 2$.
Denote by $B_1 \In \R^m$ the unit ball.
\begin{thm}\label{thm: smooth estimates2d}
For $m=2$. Let $N$ be a compact Riemannian manifold and $(\phi,\psi)$ a weakly Dirac-harmonic map from the unit disc $B_1\In \R^2$ into $N$.
For any $k\in \mathbb{N}$ there exist
$\epsilon = \epsilon(N)>0$ and $C=C(k,N)>0$ such that if
$\|\D \phi\|_{L^{2}(B_1)} \leq \epsilon,$
then
\bee \label{smooth-interior-estimates}
\|\D^k \phi\|_{L^{\infty}(B_{\fr12})} + \|\D^k \psi\|_{L^{\infty}(B_{\fr12})} \leq C\|\D \phi\|_{L^{2}(B_1)}(1 + \|\psi\|_{L^{4}(B_1)}). \eee

\end{thm}
\begin{rmk}\label{remark m=2}
In dimension $m=2$, the local smooth interior estimates as in \eqref{smooth-interior-estimates} were obtained in \cite{CJLW2} using classical methods but for smooth Dirac-harmonic maps and assuming also that $\|\psi\|_{L^4(B_1)} \leq \epsilon$.
\end{rmk}

\begin{thm}\label{thm: smooth estimates}
Let $m\geq 3$ and $N$ be a compact Riemannian manifold. For any $k\in \mathbb{N}$, there exist $\epsilon = \epsilon (m,N)>0$ and $C=C(m,N,k)>0$ such that
if $(\phi,\psi)$ is a weakly Dirac-harmonic map from the unit ball $B_1\In \R^m$ into $N$ (see Section \ref{recap} for a definition) satisfying
\bee\label{smallmorrey}  \|\D \phi\|_{M^{2,m-2}(B_1)} + \|\psi\|_{M^{4,m-2}(B_1)} \leq \epsilon
\eee
then
$$\|\D^k \phi\|_{L^{\infty}(B_{\fr12})} + \|\D^k \psi\|_{L^{\infty}(B_{\fr12})} \leq C(\|\D \phi\|_{M^{2,m-2}(B_1)} + \|\psi\|_{M^{4,m-2}(B_1)}),$$
where the Morrey norms appearing here are defined in appendix \ref{def and not}.
\end{thm}

\begin{rmk}\label{remark general domain2}
We can conclude the same results as in Theorem \ref{thm: smooth estimates} for weakly Dirac-harmonic maps from $(B_1,g_1)$ equipped with a spin structure
into $N$, where $g_1$ is some smooth Riemannian metric. This is because the crucial improvement to the regularity comes from Theorem \ref{thm: main2},
thus taking into account remark \ref{remark general domain} and equation \eqref{map1} we can essentially apply the same proof.
We also mention that this recovers the $\epsilon$-regularity estimates for weakly harmonic maps (setting $\psi\equiv 0$).
\end{rmk}

\begin{rmk}
Under the smallness assumption \eqref{smallmorrey}, the regularity for weakly Dirac-harmonic maps was proved by Wang-Xu \cite{WX} using a different method.
More precisely, they first derive the H\"older continuity for $\phi$ by combining H\'{e}lein's technique of moving
frames \cite{helein_conservation} for harmonic maps and Rivi\`ere-Struwe's Coulomb gauge theorem in the Morrey space setting \cite{riviere_struwe}, and then
obtain the $C^{1, \alpha}$-regularity for $(\phi,\psi)$ by adapting the hole-filling-type argument by Giaquinta-Hildebrandt \cite{GH} for harmonic maps.
Finally, the higher order regularity follows from the standard bootstrapping argument.
\end{rmk}

\begin{rmk}
We also point out here that it is possible to prove these estimates without the smallness assumption on $\|\psi\|_{M^{4,m-2}}$,
however since it requires more involved techniques and the result does not directly benefit us in this paper, we postpone the details to a later work.
\end{rmk}

\paragraph{Layout of the paper}
In section \ref{impPDE} we improve the regularity for suitably generalised elliptic systems appearing in \cite{Zhu} - Theorem \ref{thm: main2} -
following the ideas developed in \cite{Sh, thesis}. We also state consequences of this theorem for homogeneous Dirichlet and Neumann boundary value
problems for the same system - Theorem \ref{thm: mainPDEhom}, in order to prove Theorem \ref{thm: mainPDE}. In section \ref{recap} we recall the set-up
for Dirac-harmonic maps and the related free boundary problem developed in \cite{CJWZ}. Section \ref{highreg} contains the proofs of all of the Theorems
concerning Dirac-harmonic maps mentioned in the introduction. In the same section we consider interior and boundary regularity for systems of spinors solving
a Dirac equation, where the boundary value problem is of chirality-type - Theorems \ref{spinor_thm2} and \ref{spinor_thm1}.
In the appendix we give a brief overview of the classical function spaces we require, along with results on Hodge decompositions, the Coulomb gauge construction
of Rivi\`ere-Struwe \cite{riviere_struwe} and some classical boundary value problem estimates for the Laplace and Cauchy-Riemann operator.

{\bf Acknowledgments.}  This work was started whilst both authors were supported by The Leverhulme Trust at the University of Warwick. The second author has received funding from the European
Research Council under the European Union -- Seventh Framework Programme (FP7/2007-2013) / ERC grant agreement no. 267087. The first author has also received funding from the European Research Council STG agreement number P34897. The first author would like to thank J\"urgen Jost and the Max Plank Institute for Mathematics in the Sciences for their hospitality during the completion of this work.

\section{Supporting PDE results}\label{impPDE}

The first main PDE theorem in this section - Theorem \ref{thm: main2}, is an improved interior regularity result for a system of critical PDE.
Systems of this form were introduced and first studied by Rivi\`ere-Struwe \cite{riviere_struwe}, however in the present form the improved regularity results,
and generalisations thereof stem from works of the authors \cite{Zhu} and \cite{Sh}. The proof of Theorem \ref{thm: main2} follows precisely the lines set out
in \cite{riviere_struwe} and \cite{Zhu} in obtaining H\"older regularity and we follow \cite{Sh} in proving the higher integrability. We provide a sketch of the
proof of this theorem for the sake of clarity, referring to \cite{Sh} and \cite{thesis} for the fine detail.

\begin{thm}[]
\label{thm: main2}
Let $m, d \geq2$, $0<\Lambda <\infty$ and $\fr{m}{2}<p<m$. For any $F,G\in L^{\infty}\cap M_1^{2,m-2}(B_1, gl(d))$,
$A\in L^{\infty}\cap M_1^{2,m-2}(B_1, GL(d))$, $\Omega\in M^{2,m-2}(B_1,so(d)\otimes \wedge^1 \R^m)$,
$\zeta \in M_1^{2,m-2}(B_1,gl(d)\otimes\wedge^2 \R^m)$, $f\in L^p(B_1,\R^d)$ and $u\in M_1^{2,m-2}(B_1,\R^d)$ weakly solve
\begin{equation*}
\dv(A\ed u) = \langle \Omega , A \ed u\rangle + \langle F\dv \zeta , G \ed u\rangle + f
\end{equation*}
with
\begin{equation*}\label{eqn: Abound}
\Lambda^{-1}|\xi| \leq |A(x)\xi| \leq \Lambda |\xi| \,\,\,\,\,\text{for a.e. $x\in B_1$}, \,\,\,\text{for all $\xi\in \R^d$}
\end{equation*}
and
$
\|F\|_{L^{\infty}(B_1)} + \|G\|_{L^{\infty}(B_1)}+\|\D F\|_{M^{2,m-2}(B_1)}+\|\D G\|_{M^{2,m-2}(B_1)}  \leq \Lambda$. Then there exist $\epsilon=\epsilon(m,d,\Lambda, p)>0$ and $C=C(m,d,\Lambda, p)>0$ such that whenever $\|\Om\|_{M^{2,m-2}(B_1)}+\|\D \zeta\|_{M^{2,m-2}(B_1)}+\|\D A\|_{M^{2,m-2}(B_1)}\leq \epsilon$ then
$$\|\D^2 u\|_{M^{\fr{2p}{m},m-2}(B_{\fr12})} + \|\D u\|_{M^{\fr{2p}{m-p},m-2}(B_{\fr12})} \leq C(\| u\|_{L^1(B_1)} + \|f\|_{L^p(B_1)}).$$
\end{thm}
\begin{rmk}\label{rmkthmmain2}

 Under the conditions of Theorem \ref{thm: main2}, if $\zeta \equiv 0$ and $f\in L^r$ for all $r<m$, then we have that for all
$2< q<\infty$ there exist $\epsilon = \epsilon(m,d,\Lambda,q)>0$ and $C=C(m,d,\Lambda,q)>0$ such that whenever $\|\Om\|_{M^{2,m-2}(B_1)}+\|\D A\|_{M^{2,m-2}(B_1)}\leq \epsilon$ then
$$ \|\D^2 u\|_{M^{\fr{2q}{q+2},m-2}(B_{\fr12})} +\|\D u\|_{M^{q,m-2}(B_{\fr12})} \leq C(\| u\|_{L^1(B_1)}+\|f\|_{L^{\fr{mq}{m+q}}(B_1)}).$$
To see this notice that $p:=\fr{mq}{2+q} \in \left(\fr{m}{2},m\right)$ and apply the theorem for $p=p(q)$ where $q=\fr{2p}{m-p}.$
\end{rmk}

\begin{rmk}
Note that we could have written the PDE as
$$\dv(A\ed u) = \langle \ti{\Omega}  , A \ed u\rangle  + f$$
where $\ti{\Omega}= \Om + F\dv \zeta G A^{-1}$ is not antisymmetric. Thus it is clear that in order to have the improved regularity ``it is only necessary for the part of $\Om$ that has divergence to be anti-symmetric" as has been noted previously: \cite{Zhu} and \cite{schikorra_frames}. We also point out that there have been further generalisations of this type of theorem in different directions - the interested reader should consult the works of Roger Moser \cite{moser_regularity} and Armin Schikorra  \cite{schikorra_higher}.

\end{rmk}

In order to prove Theorem \ref{thm: main2} we first require a small generalisation of a result of Rivi\`ere-Struwe \cite{riviere_struwe} concerning
the H\"older regularity of such weak solutions:

\begin{thm}[\cite{riviere_struwe, rupflin, Zhu, schikorra_frames}]\label{thm: holder}
Let the set-up be as in Theorem \ref{thm: main2}
then there exist $\epsilon=\epsilon(m, n,\Lambda, p)>0$ and $C=C(m, n,\Lambda, p)>0$ such that whenever
$
\|\Om\|_{M^{2,m-2}(B_1)} + \| \D \zeta\|_{M^{2,m-2}(B_1)} +\|\D A\|_{M^{2,m-2}(B_1)} \leq \epsilon$
then $u\in C^{0,\gamma}_{loc}$, for $\gamma = 2-\fr{m}{p}$ with  \begin{equation*}
[u]_{C^{0,\gamma}(B_{\fr12})}\leq C (\|\D u\|_{M^{2,m-2}(B_1)} + \|f\|_{L^p(B_1)}).
\end{equation*}
\end{thm}

\begin{rmk}\label{remark general domain}
Following the proof of the above one can check that even if $f=\ed^\ast H$ for $H\in L^q(B_1, \bigwedge^1 \R^d)$ and $q>m$ we get that $u\in C^{0,\gamma}_{loc}$, for $\gamma = 1-\fr{m}{q}$ with  \begin{equation*}
[u]_{C^{0,\gamma}(B_{\fr12})}\leq C (\|\D u\|_{M^{2,m-2}(B_1)} + \|H\|_{L^q(B_1)})
\end{equation*}
as has been pointed out elsewhere \cite{JWZ}.

We point out here that both Theorem \ref{thm: main2} and Theorem \ref{thm: holder} remain true if we consider the domain to be equipped with a Lipschitz Riemannian metric $(B_1 , g_1)$ where $g_1$ is arbitrary. Of course in this case $\ed^{\ast}_{g_1}$ would depend on $g_1$ as would the inner product $\langle, \rangle_{g_1}$. We will not give a proof of this fact here as the proof follows the same lines as the Euclidean setting with some added technicalities.
\end{rmk}
\begin{rmk}
The conclusion of Theorem \ref{thm: main2} is stronger than that of Theorem \ref{thm: holder} since we have (see \cite[Chapter III, Theorem 1.2]{giaquinta}),
$[u]_{C^{0,\gamma}}\leq C\|\D u\|_{M^{\fr{2p}{m-p},m-2}}. $
\end{rmk}

Our second main PDE result in this section is to present estimates up to the boundary when we have zero Dirichlet and/or Neumann conditions.
Denote by $B_1^+:=B_1\cap \{x_m \geq 0\}$ the upper half unit ball and set $I:= \{(x_1,...,x_{m-1},0) \in \R^m |  -1 < x_1,..., x_{m-1} < 1 \}$.

\begin{thm}\label{thm: mainPDEhom}
Let $m, d \geq2$, $0\leq s\leq d$, $0<\Lambda <\infty$ and $\fr{m}{2}<p<m$. For any $A\in L^{\infty}\cap M_1^{2,m-2}(B^+_1, GL(d))$,
$\Omega\in M^{2,m-2}(B^+_1,so(d)\otimes \wedge^1 \R^m)$, $f\in L^p(B^+_1,\R^d)$ and for any $u\in M_1^{2,m-2}(B^+_1,\R^d)$ weakly solving
\begin{eqnarray}
\dv(A\ed u) &=& \langle \Omega , A \ed u\rangle  + f \quad \text{in}\quad B_1^+,   \nn \\
\pl{u^i}{\overrightarrow{n}} &=& 0  \quad \text{on $I$}\quad 1\leq i \leq s,     \label{eqn: bchom}\\
u^j &=& 0 \quad \text{on $I$}\quad s+1 \leq j \leq d,
\end{eqnarray}
with
\begin{equation*}
\Lambda^{-1}|\xi| \leq |A(x)\xi| \leq \Lambda |\xi| \,\,\,\,\,\text{for a.e. $x\in B_1^+$}, \,\,\,\text{for all $\xi\in \R^d$}
\end{equation*}
and satisfying the compatibility condition that $A(x)$ commutes with $R$ for almost every $x\in B_1^+$ where $R$ is as in \eqref{R},
there exist $\epsilon=\epsilon(m,d,\Lambda, p)>0$ and $C=C(m,d,\Lambda, p)>0$ such that whenever
$
\|\Om\|_{M^{2,m-2}(B^+_1)}+
\|\D A\|_{M^{2,m-2}(B^+_1)}\leq \epsilon$ then
$$\|\D^2 u\|_{M^{\fr{2p}{m},m-2}(B^+_{\fr12})} + \|\D u\|_{M^{\fr{2p}{m-p},m-2}(B^+_{\fr12})} \leq C(\| u\|_{L^1(B^+_1)} + \|f\|_{L^p(B^+_1)}).$$

\end{thm}
\begin{rmk}
There are versions of this theorem with inhomogeneous boundary conditions in all dimensions, i.e. analogues of Theorem \ref{thm: mainPDE} when $m\geq 3$, however we will leave the precise formulations to a later work as they require more involved techniques and do not benefit us here.
\end{rmk}

\begin{rmk}
As above we can make sense of this theorem when $B_1^+$ is equipped with a Riemannian metric $g_1^+$, but we require that it is in block form:
$g_1^+ = \twomat{\{(g_1)_{ij}^+\}_{1\leq i,j\leq m-1}}{0}{0}{(g_1)^+_{mm}}.$
The reader should think of this condition as taking local Fermi coordinates with respect to $\partial M\In M$ about some boundary point $p$
which we can always do, i.e.,  about $p\in \partial M$ we can find a chart $\phi : U\In M\to B_1^+$ such that $\phi(U\cap \partial M) = I$ and
the metric is in the form above.
\end{rmk}

In dimension $m=2$, we are able to deal with the case of inhomogeneous Neumann and Dirichlet boundary constraints - our main PDE result in the paper Theorem \ref{thm: mainPDE}.

\subsection{Proof of Theorems  \ref{thm: mainPDEhom} and \ref{thm: mainPDE}}

{\bf Proof of Theorem \ref{thm: mainPDEhom}:}
First we define reflections in the domain and target: In the domain we denote $\rho:\R^m\to \R^m$ to be
$\rho(x):=(x_1,\dots, x_{m-1},-x_m)$ and in the target, $\sigma :\R^d\to \R^d$ defined by $\sigma(y)=(y^1,\dots, y^s, -y^{s+1},\dots,-y^d)$ -
note that $R$ is the matrix describing both $\sigma$ and $D \sigma$. Now we extend the objects appearing in the Theorem as follows:
$$u_{E}(x):=\twopartdef{u(x)}{x\in B_1^+}{Ru(\rho(x))}{x\in B_1^-}\in W^{1,2}(B_1,\R^d),$$
$$A_E(x):=\twopartdef{A(x)}{x\in B_1^+}{A(\rho(x))}{x\in B_1^-}\in W^{1,2}(B_1,GL(d)),$$
$$\Om_E(x):=\twopartdef{\Om(x)}{x\in B_1^+}{R(\rho^{\ast}\Om)(x)R}{x\in B_1^-}\in L^2(B_1,so(d)\otimes\wedge^1 \R^m),$$
$$f_E(x):=\twopartdef{f(x)}{x\in B_1^+}{Rf(\rho(x))}{x\in B_1^-}\in L^p(B_1,\R^d).$$

Now, for any $\eta\in C_c^{\infty}(B_1,\R^d)$ we decompose it as follows
$$\eta(x) = \eta_1(x) + \eta_2(x) = \fr12(\eta(x) +R\eta (\rho(x))) + \fr12(\eta(x) - R\eta(\rho(x)))$$ where we have
$R\eta_1(\rho(x)) = \eta_1(x)$ and $R\eta_2(\rho(x)) = -\eta_2 (x).$
As a preliminary calculation we check that
\begin{eqnarray}
\int_{B_1^+} (\la \Om(x),A(x)\ed u(x)\ra + f(x),\eta_1 (x)) &=& \int_{B_1^-} (\la \rho^{\ast}\Om (x), A(\rho(x))\ed (u\circ\rho)(x)\ra + f(\rho(x)), \eta_1(\rho(x))) \nn\\
&=&  \int_{B_1^-} (\la \rho^{\ast}\Om (x), A(\rho(x))\ed (u\circ\rho)(x)\ra + f(\rho(x)), R\eta_1(x)) \nn\\
&=& \int_{B_1^-} (\la \Om_E (x), A_E(x)\ed u_E(x)\ra + f_E(x),\eta_1(x)). \nn
\end{eqnarray}
Therefore by \eqref {eqn: bchom} and the fact that $\eta_1^j = 0$ on $I$ for $s+1 \leq j\leq d$
we have
\bee
\int_{B_1} \left( A_E\ed u_E ,\ed \eta_1\right)&=&  \int_{B_1^+}\left( A \ed u , \ed \eta_1\right) + \int_{B_1^-} \left (A(\rho(x)) R\ed (u\circ \rho) (x), \ed \eta_1 (x)\right)\nn\\
&=&  \int_{B_1^+} \left( \la \Om, A\ed u\ra +f,\eta_1 \right) + \int_{B_1^-}  \left(A(\rho(x)) R\ed (u\circ \rho) (x), R\ed (\eta_1\circ \rho) (x)\right)\nn\\
&=&  \int_{B_1^+} \left( \la \Om, A\ed u\ra +f,\eta_1 \right) + \int_{B_1^+}  \left(A(x)\ed u(x),  \ed \eta_1 (x)\right) \nn\\
&=& 2\int_{B_1^+} \left( \la \Om, A\ed u\ra +f,\eta_1 \right)\nn\\
&=& \int_{B_1}(\la \Om_E (x), A_E(x)\ed u_E(x)\ra + f_E(x),\eta_1(x)).\nn
\eee
We also check that
\bee
\int_{B_1} \left( A_E\ed u_E ,\ed \eta_2\right)&=&  \int_{B_1^+}\left( A \ed u , \ed \eta_2\right) + \int_{B_1^-} \left (A(\rho(x)) R\ed (u\circ \rho) (x), \ed \eta_2 (x)\right)\nn\\
&=& \int_{B_1^+}\left( A \ed u , \ed \eta_2\right) - \int_{B_1^-} \left (A(\rho(x)) R\ed (u\circ \rho) (x), R\ed (\eta_2\circ\rho) (x)\right)\nn\\
&=& \int_{B_1^+}\left( A \ed u , \ed \eta_2\right) - \int_{B_1^+}\left( A(x) \ed u (x) , \ed \eta_2 (x)\right) \nn\\
&=& 0.\nn
\eee
On the other hand, a simple calculation gives
\begin{eqnarray}
\int_{B_1^+} (\la \Om(x),A(x)\ed u(x)\ra + f(x),\eta_2 (x)) &=& \int_{B_1^-} (\la \rho^{\ast}\Om (x), A(\rho(x))\ed (u\circ\rho)(x)\ra + f(\rho(x)), \eta_2(\rho(x))) \nn\\
&=&  - \int_{B_1^-} (\la \rho^{\ast}\Om (x), A(\rho(x))\ed (u\circ\rho)(x)\ra + f(\rho(x)), R\eta_2(x)) \nn\\
&=&  - \int_{B_1^-} (\la \Om_E (x), A_E(x)\ed u_E(x)\ra + f_E(x),\eta_2(x)). \nn
\end{eqnarray}
It follows that
\bee \int_{B_1}(\la \Om_E (x), A_E(x)\ed u_E(x)\ra + f_E(x),\eta_2(x))=0. \nn
\eee
Therefore we have proved that for all such $\eta$ we have
\bee
\int_{B_1} \left( A_E\ed u_E ,\ed \eta\right)=\int_{B_1}(\la \Om_E (x), A_E(x)\ed u_E(x)\ra + f_E(x),\eta(x)). \nn
\eee
By applying Theorem \ref{thm: main2} we are done.         \hfill{$\square$}

Before we begin the proof of Theorem \ref{thm: mainPDE} we first recall a fundamental result of Agmon-Douglis-Nirenberg \cite{agmon_doug_niren},
the proof of which we leave to the reader - see \cite[Theorem 3.4]{wehrheim}. We point out here that the below result holds in all dimensions but we only require the result on $B_1^+\In \R^2$.
\begin{thm} \label{Gn}
Let $G\in W^{1,p}(B_1^+)$ for $1<p<\infty$. Then there exist $C=C(p)>0$ and $\hat{u}_G \in W^{2,p}(B_1^+)$ such that
$\pl{\hat{u}_G}{\overrightarrow{n}} = G|_{I}$ and $\|\hat{u}_G\|_{W^{2,p}}\leq C\|G\|_{W^{1,p}}.$
\end{thm}

{\bf Proof of Theorem \ref{thm: mainPDE}:}
By the definition of the spaces $W_{\partial}^{k,p}$ and the above theorem we can find $\hat{u}\in W^{2,p}(B_1^+)$ such that
\begin{eqnarray*}
&\pl{\hat{u}^i}{\overrightarrow{n}}=\frak{g}^i& \quad \text{on $I$}\quad 1\leq i \leq s\\
&\hat{u}^j = \frak{k}^{j}& \quad \text{on $I$}\quad s+1 \leq j \leq d,
\end{eqnarray*} and $C=C(p)>0$ with
$\|\hat{u}\|_{W^{2,p}}\leq C(\|\frak{g}\|_{W_{\partial}^{1,p}}+\|\frak{k}\|_{W_{\partial}^{2,p}}).$
Defining $v:=u-\hat{u}$ and $F:=f-\dv(A\ed \hat{u}) + \la \Om, A\ed \hat{u}\ra$ the reader can check that we have
$$\dv (A\ed v) = \la \Om, A\ed v\ra + F \quad \text{on $B_1^+$}$$ with $v$ satisfying the boundary conditions \eqref{eqn: bchom} and
$$\|F\|_{L^p(B_1^+)} \leq C(\|f\|_{L^p(B_1^+)} + \|\frak{g}\|_{W_{\partial}^{1,p}}+\|\frak{k}\|_{W_{\partial}^{2,p}}).$$
Now we apply Theorem \ref{thm: mainPDEhom} to $v$ for $m=2$ which eventually gives the desired estimate for $u$.          \hfill{$\square$}

\subsection{Proof of Theorems \ref{thm: holder} and \ref{thm: main2}}
{\bf Proof of Theorem \ref{thm: holder}:}
The proof is an easy application of the ideas set out in \cite{riviere_struwe, rupflin}, see also \cite{schikorra_frames, Zhu}.
Setting $T:=P^{-1}A$ where $P$ is given by the Morrey version of the Coulomb gauge \cite[Lemma 3.1]{riviere_struwe} we end up with the system
$$\dv (T \ed u) = \langle\dv \eta ,T \ed u\rangle + \langle P^{-1}F \dv \zeta , GA \ed u\rangle + P^{-1}f$$
and$$\ed (T\ed u) = \ed T\wedge \ed u.$$

We have $\|T\|_{L^{\infty}} \leq \Lambda,$ $\fr1{\Lambda} |\ed u| \leq |T\ed u| \leq \Lambda |\ed u|$ and $\|\D T\|_{M^{2,m-2}(B_1)} \leq C\epsilon.$ Also $\eta = \ast \xi$ where $\xi$ is as in Lemma 3.1 of \cite{riviere_struwe}. Hence $\eta \in W_0^{1,2}(B_1, so(n) \otimes \wedge^2 \R^m)$ and $\ed \eta = 0.$ Moreover there exists $C=C(m,n)>0$ such that $\|\D \eta\|_{\M^{2,m-2}(B_1)} \leq C\epsilon.$ The proof now follows the steps precisely as they are laid out in the proof of Theorem 1.1 from \cite{riviere_struwe}, see also Proposition 2.1 in \cite{rupflin}, or Theorem 1.2 in \cite{Zhu}.
\hfill $\square$

{\bf Proof of Theorem \ref{thm: main2}:}
We first assume the following:
\begin{prop}\label{prop: improved decay}
Let the setup be as in Theorem \ref{thm: main2}. Then there exits $\epsilon = \epsilon (d,m,p,\Lambda)>0$ such that whenever $\|\Om\|_{M^{2,m-2}(B_1)}+\|\D \zeta\|_{M^{2,m-2}(B_1)} +\|\D A\|_{M^{2,m-2}(B_1)}\leq \epsilon$ then $\D u \in M_{loc}^{2,m-2+ 2\gamma}(B_1,\R^d)$.

\end{prop}

We now follow precisely the bootstrapping argument of \cite{Sh} (see also \cite[Chapter 6]{thesis}) therefore we do the first few lines and then refer the reader there for the full proof.
Setting $\Theta: = A^{-1}\ed A + A^{-1}\Om A + A^{-1}F\dv \zeta G \in M^{2,m-2}(B_1, gl(d)\otimes \wedge^1 \R^m)$ we have that $\|\Theta\|_{M^{2,m-2}(B_1)} \leq C(\Lambda)\epsilon$ and $u$ is a solution to:
$$-\Dl u = \langle\Theta , \ed u\rangle + f.$$ Using the fact that now we know $\D u$ has better than $L^2$ integrability we can simply bootstrap this information back into our PDE.
This completes the proof of Theorem \ref{thm: main2}.   \hfill$\square$

{\bf Proof of Proposition \ref{prop: improved decay}:}
We can of course apply Theorem \ref{thm: holder} to our situation in Theorem \ref{thm: main2}. Now we assume the following:
\begin{prop}\label{prop: decay1}
With the set-up as in Theorem \ref{thm: main2}; let $\dl>0$, then there exist $\epsilon=\epsilon(d,m,p, \Lambda)>0$ small enough and $C=C(\dl,m,d,p,\Lambda)>0$ such that when
$\|\Om\|_{M^{2,m-2}(B_1)}+\|\D \zeta\|_{M^{2,m-2}(B_1)}+\|\D A\|_{M^{2,m-2}(B_1)}\leq \epsilon$
we have the following estimate ($\gamma = 2-\fr{m}{p}$)
\begin{equation*}
\|A\D u\|_{L^2(B_r)}^2 \leq C(\dl)(\epsilon^2[u]_{C^{0,\gamma}(B_1)}^2 +\|f\|_{L^p(B_1)}^2)+ r^m(1+\dl)\|A\D u\|_{L^2(B_1)}^2.
\end{equation*}

\end{prop}

Assuming Proposition \ref{prop: decay1} we follow the argument for the proof of  \cite[Lemma 7.3]{Sh_To}, again the full details are given in \cite{Sh} or
 \cite[Chapter 6]{thesis}, however one should follow the argument by replacing $\D u$ with $A\D u$. This finishes the proof of Proposition \ref{prop: improved decay}. \hfill$\square$

{\bf Proof of Proposition \ref{prop: decay1}:}
This section of the proof also follows exactly the ideas set out in \cite{Sh} and \cite[Chapter 6]{thesis} but we include some of the details since there are a couple of added technicalities here: We will use the Coulomb gauge \cite[Lemma 3.1]{riviere_struwe} in order to re-write our equation; so setting $\epsilon$ small enough we have
\begin{equation*}
\dv (P^{-1}A\ed u) = \langle\dv \eta , P^{-1}A\ed u\rangle + \langle P^{-1}F\dv \zeta,  G\ed u\rangle +P^{-1}f
\end{equation*}
and
\begin{equation*}
\ed (P^{-1}A\ed u) = (\ed (P^{-1}A) \wedge \ed u).
\end{equation*}
Letting $S=P^{-1}A$ we have
\begin{equation*}
\dv (S\ed u) = \langle\dv \eta , S\ed u\rangle +\langle P^{-1}F\dv \zeta, G\ed u\rangle+P^{-1}f
\end{equation*}
and
\begin{equation*}
\ed (S\ed u) = (\ed S \wedge \ed u)
\end{equation*}
where $\frac{1}{\Lambda}|\xi|\leq |S(x)\xi|\leq \Lambda |\xi|$
for all $\xi \in \R^n$ and almost every $x\in B_1$. We also have
$\|\D S\|_{M^{2,m-2}(B_1)} \leq \epsilon (\Lambda +1).$ As before; $\eta \in W_0^{1,2}(B_1, so(n) \otimes \wedge^2 \R^m)$ and $\ed \eta = 0.$ Moreover there exists $C=C(m,n)>0$ such that $\|\D \eta\|_{\M^{2,m-2}(B_1)} \leq C\epsilon.$

We can also set $\epsilon$ small enough in order to apply Theorem \ref{thm: holder} so that $u \in C^{0,\gamma}$ where $\gamma = 2-\fr{m}{p}\in (0,1)$.
 Now we wish to extend the quantities arising above in the appropriate way: First of all we may extend $\eta$ by zero. We also extend $S-S_{1,\underline{0}}$
to $\tilde{S} \in W^{1,2}\cap L^{\infty}(\R^m)$ and finally $u$ to $\tilde{u} \in C^{0,\gamma}(\R^m)$ where each has compact support in $B_2$
(we may assume $u \in C^{0,\gamma}(\overline{B}_1)$).

Note that we have $\|\D \tilde{S}\|_{L^2} \leq C\|\D S\|_{L^2(B_1)} \leq C\epsilon (\Lambda +1)$ by Poincar\'e's inequality and $\D \tilde{S} = \D S$ in $B_1$.
 We also have $\tilde{u} \in C^{0,\gamma}(\R^m)$ with $\|\tilde{u}\|_{C^{0,\gamma}} \leq C \|u\|_{C^{0,\gamma}}$ and (since we may assume $\int u = 0 $)
we have $\|\tilde{u}\|_{C^{0,\gamma}} \leq C [u]_{C^{0,\gamma}}$, moreover $\tilde{u} = u$ in $B_1$. All the constants here come from standard extension
operators and are independent of the function, see for instance \cite{gt}. We also extend $\zeta -\zeta_{1,\underline{0}}$ to $\ti{\zeta}\in M_1^{2,m-2}(\R^m)$
 with $\D \ti{\zeta}=\D\zeta$ in $B_1$ and
$\|\D \ti{\zeta}\|_{M^{2,m-2}}\leq C\|\D \zeta\|_{M^{2,m-2}(B_1)} \leq C\epsilon$. (see appendix \ref{def and not}).
We also extend $P,F,G$ using appendix \ref{def and not}.

Now we use the $L^2$ Hodge decomposition in order to write $S\ed u = \ed a + \dv b + h$ with $a\in W_0^{1,2}(B_1, \R^n)$, $b\in W_N^{1,2}(B_1 \wedge^2 \R^m \otimes \R^n)$ (the normal part of $b$ vanishes at the boundary) has $\ed b = 0$ and $h$ is a harmonic one form with $$\|\ed a\|_{L^2(B_1)}^2 + \|\dv b\|_{L^2(B_1)}^2 + \|h\|_{L^2(B_1)}^2 = \|S \ed u\|_{L^2(B_1)}^2.$$

Notice that we have $$\Dl a = \langle\dv \eta ,S\ed u\rangle +\langle P^{-1}F\dv \zeta, G\ed u\rangle+P^{-1}f
$$ and $\Dl b = \ed S \wedge \ed u$ weakly. We proceed to estimate $S\D u \in L^2$ by estimating $\|\ed a\|_{L^2}$, $\|\dv b\|_{L^2}$ and using standard properties of harmonic functions in order to deal with $\|h\|_{L^2}$. The next part of the proof follows from \cite{Sh} or \cite[Chapter 6]{thesis} with a few extra terms involving our extensions above, however we still obtain:
\begin{equation}\label{esta}
\|\D a\|_{L^2(B_1)} \leq C (\epsilon[u]_{C^{0,\gamma}(B_1)} + \|f\|_{L^p(B_1)}).
\end{equation}
and
\begin{equation}\label{b}
\|\dv b\|_{L^2(B_1)} \leq C\epsilon[u]_{C^{0,\gamma}(B_1)}.
\end{equation}

We now use the fact that $h$ is harmonic giving that the quantity $r^{-m}\|h\|_{L^2(B_r)}^2 $ is increasing to give
\begin{eqnarray*}
\|h\|_{L^2(B_r)}^2 \leq r^m \|h\|_{L^2(B_1)}^2 \leq r^m\|S\ed u\|_{L^2(B_1)}^2 = r^m\|A\ed u\|_{L^2(B_1)}^2
\end{eqnarray*}
where the last equality follows because $P$ is orthogonal.

Going back to our original Hodge decomposition we see that (using Young's inequality, the orthogonality of $P$, \eqref{esta} and \eqref{b})
\begin{eqnarray*}
\|A\ed u \|_{L^2(B_r)}^2 &=& \|S\ed u \|_{L^2(B_r)}^2 \\
&\leq& (\|h\|_{L^2(B_r)} + \|\ed a\|_{L^2(B_r)} + \|\dv b\|_{L^2(B_r)})^2 \\
&\leq& (1+\dl)\|h\|_{L^2(B_r)}^2 + C({\dl}) (\|\ed a\|_{L^2(B_r)} + \|\dv b\|_{L^2(B_r)})^2 \\
&\leq& (1+\dl)r^m\|A\ed u\|_{L^2(B_1)}^2 + C(\dl)(\epsilon^2[u]_{C^{0,\gamma}(B_1)}^2 +\|f\|_{L^p(B_1)}^2).
\end{eqnarray*}
This completes the proof. \hfill$\square$

\section{Dirac-harmonic maps and the free boundary value problem}\label{recap}

\subsection{Dirac-harmonic maps}\label{recapdh}
In this subsection, we shall first recall the geometric set up of Dirac-harmonic maps \cite{CJLW1, CJLW2, Z, CJWZ, WX} in general dimensions $m\geq 2$ and then
derive the Dirac-harmonic map system in some appropriate forms so that the PDE results in Section 2 can be applied.

Let $(M,g_M)$ be a Riemannian manifold of dimension $m \geq 2$ and with a fixed spin structure. Let $\Sigma M$ be the spinor bundle over $M$ with
a hermitian metric $\langle\cdot,\cdot\rangle_{\Sigma M}$ and a compatible spin connection $\nabla$. For $X\in \Gamma(TM)$ and $\psi \in \Gamma(\Sigma M)$,
denote by $X \cdot \psi$ the Clifford multiplication, satisfying:
\bee \label{clifford}  \la X \cdot \psi, \varphi \ra = - \la \psi,  X \cdot \varphi \ra, \quad    X \cdot Y \cdot \psi + Y \cdot X \cdot \psi = -2 g_M(X,Y) \psi,
\eee
for any $ X, Y \in \Gamma(TM)$ and  any $\psi, \varphi \in \Gamma(\Sigma M) $. The usual Dirac operator $\Par$ is defined by
$\Par \psi:= \gamma_{\alpha} \cdot \nabla _{\gamma_{\alpha}} \psi$ for a local orthonormal frame $\{\gamma_{\a}\}_{\a=1}^{m}$ on $M$ and
a usual spinor $\psi\in\Gamma(\Sigma M)$. The summation convention will be used throughout the paper. We refer to \cite{LM} for more spin geometric materials.

Let $(N,g)$ be a compact Riemannian manifold of dimension $d \ge 2$ and let its metric in local coordinates $\{y^i\}_{i=1}^d$ be given by $g_{ij}$,
with Christoffel symbols $\Gamma^{i}_{jk}$ and Riemannian curvature tensor $\mathcal{R}^{m}_{lij}$. Let $\phi$ be a smooth map from $M$ to $N$
and $\phi^{-1}TN$ the pull-back bundle of $TN$ under $\phi$. Consider the twisted bundle $\Sigma M\otimes \phi^{-1}TN$ with the induced metric
$\langle\cdot,\cdot\rangle_{\Sigma M\otimes \phi^{-1}TN}$ and the induced connection $\widetilde{\nabla}:=\nabla\otimes 1 + 1 \otimes \nabla^{\phi^{-1}TN}$.
In local coordinates, the section $\psi$ of $\Sigma M\otimes \phi^{-1} TN$ can be expressed by $\psi(x)= \sum_{j=1}^d\psi^j (x)\otimes \partial_j(\phi(x)),$
where $\psi^j$ is a usual spinor and $\{\partial_j = \partial y^j\}$ is the natural local basis on $N$. $\widetilde{\nabla}$ can be expressed by
$$\widetilde{\nabla} \psi= \sum_{i=1}^d\nabla \psi^i(x)\otimes \partial_i(\phi(x))
+ \sum_{i,j,k=1}^d\Gamma^i_{jk}(\phi(x)) \phi^j_{\ast}(\gamma_\a)\gamma_\a \cdot
\psi^k(x)\otimes \partial_i (\phi(x)).$$
The {\it Dirac operator along the map $\phi$} is defined by
 \begin{eqnarray*}
\Dir\psi := \g_\a\cdot \widetilde{\nabla}_{\g_\a}
\psi\nonumber = \sum_{i} \Par  \psi^i(x) \otimes \partial_i(\phi(x))
+ \sum_{i,j,k=1}^d\Gamma^i_{jk}(\phi(x)) \phi^j_{\ast}(\gamma_\a)\gamma_\a\cdot \psi^k(x) \otimes \partial_i (\phi(x)).
\end{eqnarray*}
Set
\[\mathcal{X}(M,N)  := \left \{ (\phi,\psi)\,|\, \phi \in C^\infty(M,N) \hbox{ and }
\psi \in C^\infty (\Sigma M\otimes \phi^{-1}TN) \right \}.\]
Consider the following functional on $\mathcal{X}(M,N)$:
\begin{eqnarray*}
 L(\phi,\psi)&:=&\int_M\left (|d\phi|^2+\la\psi,\Dir\psi\ra\right) \ed V_{g_{M}}.
\end{eqnarray*}
In terms of local coordinates, the corresponding Euler-Lagrange equations are:
\bee
\label{2.1} - \Delta_{g_M} \phi^{i}+\Gamma^{i}_{jk}(\phi)\la\ed \phi^{j},\ed \phi^{k}\ra_{g_{M}}
- \frac{1}{2} \mathcal{R}^{i}_{lkj}(\phi)\langle\psi^{k},\nabla\phi^{l}\cdot\psi^{j}\rangle =0,
\quad i=1,2,\ldots,d,
\eee
\bee
\label{2.2} \Par \psi^i+\Gamma^i_{jk}(\phi)\D\phi^j\cdot\psi^k=0,  \quad i=1,2,\ldots,d.
\eee
Here $\Delta_{g_M}$ is the Laplace-Beltrami operator with respect to $g_M$, $\nabla \phi^{l}:=\D^{g_M} \phi^l = \phi^{l}_{\ast}(\gamma_{\al}) \g_{\alpha}=
\phi_{\alpha}^{l} \g_{\alpha}$ and $``\cdot"$ denotes the Clifford multiplication. Solutions $(\phi,\psi) \in \mathcal{X}(M,N) $ to \eqref{2.1}, \eqref{2.2} are called Dirac-harmonic maps from $M$ to $N$.

We embed $N$ isometrically into some Euclidean space $\R^n$ via the Nash embedding theorem.
Let $A(\cdot,\cdot)$ be the second fundamental form of $N$ in $\mathbb{R}^{n}$, and $P$ the shape operator of $N$ in $\mathbb{R}^{n}$.
Set
\bee
\mathcal{A}(d\phi(\g_{\alpha}),\g_{\alpha}\cdot\psi)&=&\phi^{i}_{\alpha}\g_{\alpha}\cdot\psi^{j}\otimes
A(\partial_{y^{i}},\partial_{y^{j}}),\nn   \\
\mathcal{P}(\mathcal{A}(d\phi(\g_{\alpha}),\g_{\alpha}\cdot\psi);\psi)\nn
&=&P(A(\partial_{y^{l}},\partial_{y^{j}});\partial_{y^{i}})\langle\psi^{i},\g_{\alpha}\cdot\psi^{j}\rangle\phi_{\alpha}^{l}.
\eee
By the Gauss equation, the equations \eqref{2.1} and \eqref{2.2} can be written as (see \cite{Z, CJWZ})
\begin{equation}
\label{2.8}- \Delta_{g_M} \phi=A(\phi_{\al}, \phi_{\a})+ {\rm Re}\
\mathcal{P}(\mathcal{A}(d\phi(\g_\alpha),\g_\alpha\cdot\psi);\psi),
\end{equation}
\begin{equation}
\label{2.9}
\partial \hskip -1.8mm \slash
 \psi = \mathcal{A}(d\phi(\g_\a), \g_\a\cdot \psi ).
\end{equation}
Denote
$$H^1(M,N):=\left \{\phi\in H^1( M, \RR^K)  \   |   \   \phi (x)\in N \
a. e. \ x\in M  \right \},$$
$$\mathscr{W}^{1,\fr43}(\Sigma M\otimes \phi^{-1}TN):=   \left \{\psi\in \Gamma (\Sigma M\otimes
\phi^{-1}TN) \  | \  \int_M|\nabla \psi|^{\fr43}<+\infty, \int_M|
\psi|^{4}<+\infty \right \}.$$ Here, $\psi \in \Gamma (\Sigma
M\otimes \phi^{-1}TN)$, the spinor field along the map $\phi$,
should be understood as an $n$-tuple of spinors
$(\psi^{1},\psi^{2},...,\psi^{n})$ satisfying $ \sum_{i}\nu^{i}\psi^{i}=0, \ {\rm for\ any\ normal\
vector}\ \nu=\sum_{i}\nu^{i}\partial z^{i}\ {\rm at}\ \phi(x)$, where $\{z^i, i=1,2,...,n\}$ are the canonical coordinates of
$\mathbb{R}^{n}$. In the case that $\dim{M}=2$ we have $\mathscr{W}^{1,\fr43} = W^{1,\fr43}$ by Sobolev embedding. It is easy to verify that $L(\phi,\psi)$ is well defined for $(\phi,\psi)\in H^1(M,N)\times \mathscr{W}^{1, \fr43}(\Sigma M\otimes \phi^{-1}TN)$ (see \cite[Remark 1.3]{WX}).
Denote $$\mathcal{X}^{1,2}_{1,\fr43}(M,N) := \left \{ (\phi,\psi) \in H^1(M,N)\times \mathscr{W}^{1, \fr43}(\Sigma M\otimes
\phi^{-1}TN) \right \}$$
and extend the functional $L(\cdot, \cdot)$ to the space $\mathcal{X}^{1,2}_{1,\fr43}(M,N)$. Critical points $(\phi,\psi)\in \mathcal{X}^{1,2}_{1,\fr43}(M,N)$ of $L(\cdot, \cdot)$
are called weakly Dirac-harmonic maps from $M$ to $N$ (see \cite{CJLW2, WX}).

It is worth remarking that the equation \eqref{2.8} can be written as an elliptic system with an $L^2$-antisymmetric structure and hence the Coulomb gauge construction
in \cite{riviere_struwe} can be applied to prove the interior continuity of $\phi$ for any weakly Dirac-harmonic map $(\phi,\psi)$ from a spin Riemann surface
$M$ into a compact Riemannian manifold $N$ (see Theorem 2.1 in \cite{CJWZ} and Theorem 1.5 in \cite{WX}), extending the case of $N=S^d$
in \cite{CJLW2} and the case that $N$ is a compact hypersurface in $\R^{d+1}$ in \cite{Z}.

To see this, we consider the case of a domain $(B_1\In \R^m, g_1)$ and apply a similar procedure as in the case of a Euclidean disc done in \cite{Z, CJWZ}.
Take a local orthonormal basis $\{\g_{\a}, \a=1,...,m\}$, its dual basis $\{\theta^{\a}, \a=1,...,m\}$ and the canonical coordinates $(z^{1},z^{2},..., z^{n})$
of $\mathbb{R}^{n}$. Let $\nu_{L}$, $L=d+1,...,n$ be an orthonormal frame field for the normal bundle $T^{\bot}N$ to $N$. Denote by $\nu_{L}$ the
 corresponding unit normal vector field along the map $\phi$.  Write $\phi=\phi^{i}\partial z^{i}$, $\psi=\psi^{j}\otimes \partial z^{j}$ and
denote $\phi_{\a}=\phi_{\ast}(\g_{\a})=\phi_{x_{\a}}, \a= 1,,...,m$. Then, similarly to the case of a Euclidean disc considered in \cite{Z, CJWZ}, the equations \eqref{2.8} and \eqref{2.9}  can be written in the following
extrinsic form in terms of the orthonormal frame field $\nu_{L}$, $l=d+1,...,n$, for $T^{\bot}N$
\bee \label{map}  \Delta_{g_1}\phi ^{i}&=& \sum_{k,j,L,\a}
\left(\phi_{\alpha}^{j} \frac{\partial
\nu_{L}^{k}}{\partial z^{j}}\nu_{L}^{i}-\phi_{\alpha}^{j}
\frac{\partial \nu_{L}^{i}}{\partial z^{j}}\nu_{L}^{k}\right) \phi_{\alpha}^{k}  \nn \\
 &&+\sum_{k,j,l,L,\a} \langle\psi^{l}, \g_{\alpha}\cdot\psi^{j}\rangle
\left((\frac{\partial \nu_{L}}{\partial
z^{j}})^{\top,k}(\frac{\partial \nu_{L}}{\partial
z^{l}})^{\top,i}-(\frac{\partial \nu_{L}}{\partial
z^{l}})^{\top,k}(\frac{\partial \nu_{L}}{\partial
z^{j}})^{\top,i}\right)\phi_{\alpha}^{k},   \\
 \label{spinor}  \Par \psi ^{i} &=&  \sum_{k,j,L,\a}\frac{\partial
\nu_{L}^{k}}{\partial z^{j}} \nu_{L}^{i} \phi_{\alpha}^{k}
\g_{\alpha}\cdot\psi^{j}.
 \eee
Here $\top$ denotes the orthogonal projection $: \mathbb{R}^{n}\rightarrow T_{z}N$ and $(\cdot)^{i}$ denotes the $i$-the component of a vector
of $\mathbb{R}^{n}$. Of course, in general we cannot assume that the normal bundle of $N$ is trivial, however a simple argument using a partition of unity
allows us to consider this case without loss of generality. In particular let $\{\chi_t\}$ be a partition of unity such that the normal bundle of $N$ is trivial over the support of
each $\chi_t$ - let $\nu_{l,t}$ denote a corresponding smooth frame. Now, since the expressions above are independent of the orthonormal normal frame,
we are free to cut off any such frame using $\chi_t$ and summing up we are done - the details are left to the reader or see \cite[p.5]{moser_regularity} for the details in the case of harmonic maps.

Set
\begin{eqnarray*}
(\Om^i_k)_{\alpha}\theta^{\al}&:=& \left[\sum_{j,L}\left(\phi_{\alpha}^{j} \frac{\partial
\nu_{L}^{k}}{\partial z^{j}}\nu_{L}^{i}-\phi_{\alpha}^{j}
\frac{\partial \nu_{L}^{i}}{\partial z^{j}}\nu_{L}^{k}\right)\right]\theta^{\a} \\
&&+ \left[\sum_{l,j,L}\langle\psi^{l}, \g_{\alpha}\cdot\psi^{j}\rangle
\left((\frac{\partial \nu_{L}}{\partial
z^{j}})^{\top,k}(\frac{\partial \nu_{L}}{\partial
z^{l}})^{\top,i}-(\frac{\partial \nu_{L}}{\partial
z^{l}})^{\top,k}(\frac{\partial \nu_{L}}{\partial
z^{j}})^{\top,i}\right)\right]\theta^{\a} \\
&=& (\om^i_k)_{\al}\theta^{\al} + (F^i_j)_{\al}\theta^\al.
\end{eqnarray*}

Then we can write equation \eqref{map} in the following two equivalent forms:
\begin{equation} \label{map1}  \Delta_{g_1} \phi^{i} = \la\Omega_{k}^{i}, \ed \phi^{k}\ra_{g_1}
\end{equation}
\begin{equation}\label{map2}
 \Delta_{g_1} \phi^i = \la \om^i_k, \ed \phi^k \ra_{g_1} + f^i
\end{equation}

with $\Omega = ((\Omega^{i}_{k} )_{\alpha})_{1 \leq i,k \leq n}^{1\leq\alpha\leq m}\in L^{2}(B_1,so(n)\otimes \wedge^1 \T^{\ast}\mathbb{R}^{m}) $,
$\om =((\omega^{i}_{k} )_{\alpha})_{1 \leq i,k \leq n}^{1\leq\alpha\leq m} \in L^{2}(B_1,so(n)\otimes \wedge^1 \T^{\ast}\mathbb{R}^{m}) $ and
$f^i:=\la F^i_j, \phi^j\ra_{g_1}$ satisfying
\bee  |\Omega(x)|  \leq C (|\nabla \phi(x)| + |\psi(x)|^2),  \quad \text{for}\ a.e. \ x \in  B^+_1.  \nn\\
|\om(x) |  \leq C |\nabla \phi(x)|,  \quad \text{for}\ a.e. \ x \in  B^+_1.  \nn\\
|f(x)| \leq C |\psi(x)|^2 |\D \phi(x)|,   \quad \text{for}\ a.e. \ x \in  B^+_1.  \nn
\eee

\subsection{Free boundary value problem for Dirac-harmonic maps from surfaces}\label{recapfbvp}

In this subsection, we shall first recall the free boundary value problem for Dirac-harmonic maps from surfaces introduced in \cite{CJWZ} and then derive
the appropriate setting so that our main PDE result Theorem \ref{thm: mainPDE} can be applied.

Let $M$ be a compact spin Riemann surface with boundary $\partial M \neq \emptyset$ and let $\mathcal{S}$ be a closed $s$-dimensional submanifold of $N$. Following the notations in Section 3 of \cite{CJWZ}, we denote by $\sigma$ the geodesic reflection about $\mathcal{S}$ on a tubular neighbourhood $ \mathbf{U}_{\delta}:=\left \{y \in N  |{\rm
dist}^N(y,\cal{S})<\delta \right \}$ of $\mathcal{S}$ in $N$ and define a natural $(1,1)$ tensor $R=D \sigma$ on $\mathcal{S}$ that is compatible. In terms of the Fermi coordinates (or geodesic parallel coordinates) $\{y^i\}_{i=1}^d$ for the submanifold $\mathcal{S}$, namely,
$\{y^a\}_{a=1}^{s}$ are coordinates in $\cal{S}$ and $ \{y^\lambda\}_{\lambda=s+1}^{d}$ are the directions normal to $\cal{S}$, the tensor $R$ and
the metric $g$ on $N$ take the following forms
\bee \label{Rg}  R =  \left(
           \begin{array}{cc}
             \delta^a_b & 0\\
            0 & -\delta^\lambda_\mu \\
           \end{array}
         \right), \quad  g =  \left(
           \begin{array}{cc}
             g_{ab} & 0\\
            0 & g_{\lambda \mu} \\
           \end{array}
         \right).
\eee
Here and in the sequel, we use the following index ranges: $1\leq a,b,\cdots\leq s,  s+1\leq \lambda, \mu, \cdots\leq d, 1\leq i,j,k, \cdots\leq d.$

Let $\phi \in C^{\infty}(M,N)$ be a map from $M$ to $N$ and $\psi \in C^{\infty}(\Sigma M\otimes \phi^{-1}TN)$ a spinor field along $\phi$.
For the map $\phi$, one can impose the free boundary condition for a map in the classical sense, namely, $\phi(\partial M) \subset \mathcal{S}$.
Recall that a spinor field $\psi\in C^{\infty}(\Sigma M\otimes \phi^{-1}TN)$ along a map $\phi$ satisfies the chirality boundary conditions $\mathbf{B}_{\phi}^{\pm}$
if
\bee \label{B1} \mathbf{B}_{\phi}^{\pm}\ \psi|_{\partial M} = 0.
\eee
To reformulate the boundary condition \eqref{B1} in terms of local coordinates, we consider the case of a domain $M= \R^2_+$ equipped with the Euclidian metric $dx_1^2+dx_2^2$. Let $\gamma_{\a}= \partial x_{\a}, \a =1,2$
 be the standard orthonormal frame. A spinor field on $\R^2_+$ is simply a section $\om\in \Gamma(\Sigma \R^2_+ = \R^2_+\times \mathbb{C}^2)$ given by
$\omega = \left(
                      \begin{array}{c}
                        \omega_1 \\
                        \omega_2 \\
                      \end{array}
                    \right) :\R^2_+\to \mathbb{C}^2$ and the Clifford multiplication of $\gamma_{\a}, \a =1,2$
acting on spinors can be seen as matrix multiplication with the following identification:
\begin{equation}\label{representation}
 \gamma_1\to \begin{pmatrix} 0 & i \\ i & 0
\end{pmatrix},\qquad   \gamma_2 \to \begin{pmatrix} 0 & 1 \\
-1 & 0 \end{pmatrix}.
\end{equation}
Here, without loss of generality, we keep the representations of $\gamma_1$ and $\gamma_2$ consistent with those in \cite{CJWZ}. If exchanging $\gamma_1$
and $\gamma_2$, then their representations are consistent with those in \cite{CJLW1} and this case can be handled analogously. The usual Dirac operator $\Par$ is therefore given by
\begin{equation}\label{diracoperator}
\Par \om = \gamma_{\a} \cdot \nabla_{\gamma_{\a}} \omega =\begin{pmatrix} 0 & i \\ i & 0
\end{pmatrix}\begin{pmatrix} \frac{\partial
    \om_1}{\partial x_1} \\[2ex] \frac{\partial
    \om_2}{\partial x_1} \end{pmatrix}+\begin{pmatrix} 0 & 1 \\ -1 & 0
\end{pmatrix}\begin{pmatrix} \frac{\partial
    \om_1}{\partial x_2} \\[2ex] \frac{\partial
    \om_2}{\partial x_2} \end{pmatrix}=2i \begin{pmatrix} \frac{\partial
   \om_2}{\partial z} \\[2ex] \frac{\partial
   \om_1}{\partial \overline{z}} \end{pmatrix}.
\end{equation}
For a spinor field  (along a map $\phi$) $\psi= \psi^i \otimes \partial_i$, we write $\psi^i = \left(
                      \begin{array}{c}
                        \psi^i_+ \\
                        \psi^i_- \\
                      \end{array}
                    \right)$.
Then the chirality boundary condition \eqref{B1} becomes:
\bee  \label{localB1} \psi^i_+= \mp \ R^i_j\ \psi^j_-, \quad i
=1,2,\cdots,d,  \quad {\rm on}\ \partial M \eee
where $R=(R^i_j)$ is as in \eqref{Rg}.

Set
\bee \mathcal{X}(M,N;\mathcal{S}):=\left \{(\phi,\psi)| \phi\in
C^{\infty}(M,N), \phi(\partial M) \subset \mathcal{S};  \psi \in
C^{\infty}(\Sigma M \otimes \phi^{-1}TN), \mathbf{B}_{\phi}^{\pm}\
\psi|_{\partial M} = 0 \right \}. \nn\eee
A critical point $(\phi,\psi)$ of $L(\cdot,\cdot)$ in $\mathcal{X}(M,N;\mathcal{S})$ is called a Dirac-harmonic map from $M$ to $N$
with free boundary $\phi(\partial M)$ on $\mathcal{S}$.

Let $(\phi,\psi)$ be a Dirac-harmonic map from $M$ to $N$ with free boundary $\phi(\partial M)$ on $\mathcal{S} \subset N$.
As is done in \cite[p.1011-1012]{CJWZ}, using local coordinates $\{y^i\}_{i=1}^d$, we can take the following two types of admissible variations
$(\phi_t,\psi_{t}) \in \mathcal{X}(M,N;\mathcal{S})$:
\begin{itemize}
\item[i)]  $\phi_t\equiv \phi$ and  $\frac{d \psi_t}{dt}|_{t=0}=\xi$ is arbitrary;
\item[ii)]   $\frac{d \phi_t}{dt}|_{t=0}=\eta$ is arbitrary  and $\psi_t = \psi_t^i\otimes \partial_i(\phi_t)$ with $\psi_t^i \equiv \psi^i$
\end{itemize}
and applying direct calculations from calculus of variations to obtain the following boundary constraint:
\bee  \label{localB2} \left(2 \phi_{\overrightarrow{n}}^j - g^{jn} \langle \overrightarrow{n}\cdot\psi^l, \psi^i
\rangle \Gamma^k_{in}g_{kl}\right) \partial_j \perp T_\phi \mathcal{S},  \quad {\rm on}\ \partial M.
\eee
where $\partial_j:=\frac{\partial}{\partial y^j}$.

Furthermore, in terms of the Fermi coordinates $\{y^i\}_{i=1}^d$ about $\mathcal{S}$,  we are able to reformulate \eqref{localB2} into a simpler form.
We shall proceed as in \cite[p.1012]{CJWZ}. By the chirality boundary condition  \eqref{B1} (see also \eqref{localB1}) for $\psi$, there hold:
\begin{equation}
\label{localb}  \la \overrightarrow{n}\cdot \psi^a,\psi^b  \ra =0,
\quad \la \overrightarrow{n}\cdot \psi^\lambda,\psi^\mu  \ra =0,
\quad {\rm on}\ \partial M
\end{equation}
for $a,b=1,2,\dots,s$ and $\lambda,\mu =s+1, \dots,d.$
Since $g_{a \lambda} = 0$ (see \eqref{Rg}), one can verify that
\bee \label{christ}  g_{a b }\Gamma^a_{\lambda n}= - g_{\lambda\mu}\Gamma^\mu_{bn}.
\eee
By \eqref{clifford}, \eqref{localb} and \eqref{christ}, we have
\bee
g^{jn}g_{kl}\Gamma^k_{in}\la \overrightarrow{n}\cdot \psi^l,\psi^i\ra &=& g^{jn}g_{\lambda l}\Gamma^\lambda_{in}\la \overrightarrow{n}\cdot \psi^l,\psi^i\ra
+g^{jn}g_{a l}\Gamma^a_{in}\la \overrightarrow{n}\cdot \psi^l,\psi^i\ra   \nn\\
&=&g^{jn}g_{\lambda \mu}\Gamma^\lambda_{in}\la \overrightarrow{n}\cdot \psi^\mu,\psi^i\ra
+g^{jn}g_{a b }\Gamma^a_{in}\la \overrightarrow{n}\cdot \psi^b,\psi^i\ra    \nn\\
&=&g^{jn}g_{\lambda \mu}\Gamma^\lambda_{an}\la \overrightarrow{n}\cdot \psi^\mu,\psi^a\ra
+g^{jn}g_{a b }\Gamma^a_{\lambda n}\la \overrightarrow{n}\cdot \psi^b,\psi^\lambda\ra \nn\\
&=&g^{jn}g_{\lambda \mu}\Gamma^\lambda_{an}\la
\overrightarrow{n}\cdot \psi^\mu,\psi^a\ra+g^{jn}g_{\lambda\mu}\Gamma^\mu_{bn}\la \psi^b,\overrightarrow{n}\cdot \psi^\lambda \ra, \nn\\
&=&g^{jn}g_{\lambda \mu}\Gamma^\lambda_{an}\la
\overrightarrow{n}\cdot \psi^\mu,\psi^a\ra + g^{jn}g_{\lambda\mu}\Gamma^\lambda_{an}\la \psi^a,\overrightarrow{n}\cdot \psi^\mu \ra, \nn\\
&=&g^{jn}g_{\lambda \mu}\Gamma^\lambda_{an}\la
\overrightarrow{n}\cdot \psi^\mu,\psi^a\ra + g^{jn}g_{\lambda\mu}\Gamma^\lambda_{an} \overline{\la \overrightarrow{n}\cdot \psi^\mu,  \psi^a\ra}, \nn\\
&=&   2 g^{jn}g_{\lambda \mu}\Gamma^\lambda_{an}{\rm Re}(\la \overrightarrow{n}\cdot \psi^\mu,\psi^a\ra)     \quad \text{on $\partial M$.} \nn
\eee
Therefore, \eqref{localB2} is equivalent to
\bee \label{localB3}   (\frac{\partial \phi}{\partial \overrightarrow{n}})^\top =  g^{cd}\Gamma^\lambda_{ad}g_{\lambda\mu}{\rm Re} (\la \overrightarrow{n}\cdot \psi^\mu,\psi^a\ra)
  \partial_c \quad \text{on $\partial M$.}
 \eee
Let $P_\mathcal{S}(\cdot;\cdot)$ be the shape operator of $\mathcal{S}$ in $N$ and define
$\mathcal{P}_\mathcal{S} ( \psi^{\mu} \otimes \partial_{\mu};  \psi^{a} \otimes \partial_{a})
:= P_\mathcal{S}( \partial_{\mu};  \partial_{a})   \la \psi^{\mu},  \psi^{a} \ra $.
Then  \eqref{localB3}  is equivalent to
\begin{equation} \label{B2}
\left(\pl{\phi}{\overrightarrow{n}}\right)^{\top} = {\rm Re}(  \mathcal{P}_{\mathcal{S}}(\overrightarrow{n}\cdot \psi^{\bot};\psi^{\top})) ,
  \quad \text{on $\partial M$.}
\end{equation}
Here $\psi^{\bot}$ denotes the normal part of $\psi$ and $\psi^{\top}$ denotes the tangent part of $\psi$.

To summarise, we have the following equivalent definition:

\begin{defi} [\cite{CJWZ}, Definition 3.2] \label{bcon} $(\phi,\psi)\in \mathcal{X}(M,N;\mathcal{S})$ is
called a Dirac-harmonic map from $M$ to $N$ with free boundary $\phi(\partial M)$ on $\mathcal{S} \subset N$ if $(\phi,\psi)$ is Dirac-harmonic in $M$ -  namely they solve
\eqref{2.1}, \eqref{2.2} in $M$ - and they satisfy the boundary conditions \eqref{B1}, \eqref{B2}.
\end{defi}

Following \cite[p.1013-1014]{CJWZ}, we can use the isometric embedding $N \subset\mathbb{R}^n$ to define the free boundary conditions for
weakly Dirac-harmonic maps from surfaces. Set
\bee   \mathcal{X}^{1,2}_{1,\fr43}(M,N;\mathcal{S}) := \left \{
\ba{ll}
 (\phi,\psi)  \in H^1(M,N) \times W^{1,\fr43}(\Sigma M \otimes \phi^{-1}TN): \nn\\ \phi(x) \in \mathcal{S} \  \text{and}  \
 \mathbf{B}_{\phi}^{\pm}\ \psi (x) = 0  \text{ for  a.e. } x \in  \partial M \nn
\ea
\right \}
\eee
A critical point $(\phi,\psi)$ of $L(\cdot,\cdot)$ in $\mathcal{X}^{1,2}_{1,\fr43}(M,N; \mathcal{S})$ is
called a weakly Dirac-harmonic map with free boundary $\phi(\partial M)$ on $\mathcal{S}$.

\begin{prop}  \label{pro3.7} Let $(\phi,\psi)\in \mathcal{X}^{1,2}_{1,\fr43}(M,N;\mathcal{S})$ be a weakly Dirac-harmonic map with free boundary $\phi(\partial M)$ on $\mathcal{S}$.
Suppose that $\phi$ maps $M$ into a single coordinate neighbourhood $\{y^i\}_{i=1}^d$ of $N$ and hence our spinor can be expressed as
$\psi = \sum_{i=1}^{d} \psi^i\otimes \partial_i (\phi) $.
Then
\begin{equation} \label{weakmap}
\int_{M} \la \ed \phi, \nabla^{\phi} V\ra + \int_{M} \la\psi^k,\g_\a \cdot \psi^l\ra \la \mathcal{R}^h_{jkl} \phi_{\ast}^j(\g_{\a})\partial_h, V\ra
= \int_{\partial M}  \frac{1}{2} \la \overrightarrow{n}\cdot\psi^l, \psi^i \ra  \la  g^{jn}\Gamma^k_{in}g_{kl}   \partial_j , V\ra
\end{equation}
\begin{equation*}
\label{weakspinor}  \int_{M} \left \la \psi, \Dir \xi \right \ra = 0,  \end{equation*}
for all  $V\in H^1\cap L^{\infty}(M, \phi^{-1}TN)$ such that $V(x) \in T_{\phi(x)}\mathcal{S}$ for $a.e.\,\, x \in \partial M$ and for all
$\xi \in W^{1,\fr43}\cap L^{\infty}(\Sigma M \otimes\phi^{-1}TN)$ such that $\mathbf{B}_{\phi}^{\pm}\ \xi|_{\partial M} = 0$.

In particular, in terms of Fermi coordinates $\{y^i\}_{i=1}^d$ about $\mathcal{S}$, \eqref{weakmap} becomes
\begin{equation}\int_{M} \la \ed \phi, \nabla^{\phi} V\ra + \int_{M} \la\psi^k,\g_\a \cdot \psi^l\ra \la \mathcal{R}^h_{jkl} \phi_{\ast}^j(\g_{\a})\partial_h, V\ra
=  \int_{\partial M} ({\rm Re} \la \overrightarrow{n}\cdot \psi^\mu,\psi^a\ra ) \la  g^{cd}\Gamma^\lambda_{ad}g_{\lambda\mu}   \partial_c, V\ra \nn.
\end{equation}
\end{prop}

{\bf Proof}: Since we can use local coordinates $\{y^i\}_{i=1}^d$ on the target, similarly to the smooth case (see \cite[p.1011-1012]{CJWZ}), we can take two types of admissible variations
$(\phi_t,\psi_{t}) \in \mathcal{X}^{1,2}_{1,\fr43}(M,N;\mathcal{S})$ and apply direct calculations to complete the proof. \hfill $\square$

\begin{rmk}
The condition that $\phi$ maps into a single coordinate neighbourhood of $N$ seems rash until one considers Lemma \ref{nhdlemma}.
Notice that since $\psi\in W^{1,\fr43}$ we have that $\psi \in L^2(\partial M)$ in a trace sense, moreover since $V\in W^{1,2}\cap L^{\infty}$ when we take a trace in the $W^{1,2}$ sense (not necessarily bounded) we actually must have that the trace is bounded. Moreover the operator $\D^{\phi}$ is of course the pull back Levi-Civita connection on $\phi^{-1}TN$.
\end{rmk}

\subsection{Interior and boundary regularity for spinors}
We now state an $\epsilon$-regularity theorem for solutions to a linear Dirac system of spinors in two dimensions:
\begin{thm}\label{spinor_thm2}
Let $B_1\In \R^2$ and $\psi\in L^{4}(B_1,\mathbb{C}^2 \otimes \R^d)$, $\Gamma \in L^2(B_1,gl(d)\otimes\wedge^1 \R^2)$ solve
\begin{equation*}
\Par \psi = \Gamma_{\al} \gamma_{\al}\cdot \psi\quad \text{on $B_1$.}
\end{equation*}
Then for any $2\leq q<\infty$, $\psi\in W_{loc}^{1,\fr{2q}{2+q}}$ and there exist $\epsilon = \epsilon(q)>0$ and $C=C(q)>0$ such that whenever
$\|\Gamma\|_{L^2(B_1)} \leq \epsilon$ then
\begin{equation*}
\|\psi\|_{L^q(B_{\fr12})} + \|\D \psi\|_{L^{\fr{2q}{2+q}}(B_{\fr12})} \leq C\|\psi\|_{L^{4}(B_1)}.
\end{equation*}
\end{thm}

The above theorem has a counterpart when the spinor has a Chirality-type boundary condition;

\begin{thm}\label{spinor_thm1}
Let $B_1^+\In\R^2$ and suppose $\psi= \cvec{\psi_+}{\psi_-} \in W^{1,1}(B_1^+,\mathbb{C}^2\otimes \R^d)$ 
solves the following boundary value problem:
\begin{equation*}
\Par \psi = \Gamma_{\al} \gamma_{\al}\cdot \psi\quad\text{in $B_1^+$}
\end{equation*} and
\begin{equation*}
\mathbf{B}_R^{\pm}\psi = 0 \quad \text{on $I$}\quad \iff \quad\psi_+ = \mp R\psi_- \quad \text{on $I$}
\end{equation*}
for $\Gamma \in L^{2}(B_1^+,gl(d) \otimes \wedge^1 \R^2)$ and $R$ is as in \eqref{R}. Then for any $2 \leq q<\infty$,
$\psi\in W_{loc}^{1,\fr{2q}{2+q}}(B_1^+\cup I)$ and there exist $\epsilon=\epsilon(q)>0$ and $C=C(q)>0$ such that if
$\|\Gamma\|_{L^2(B_1^+)}\leq \epsilon$ then
\begin{equation*}
\|\D \psi\|_{L^{\fr{2q}{2+q}}(B_{\fr12}^+)} + \|\psi\|_{L^q(B_{\fr12}^+)} \leq C\|\psi\|_{L^{4}(B_1^+)}.
\end{equation*}
\end{thm}

{\bf Proof of Theorem \ref{spinor_thm2}:}
The proof of such a theorem is now standard and is essentially contained in \cite{changyouwang_dirac} so we only provide a sketch:
apply the Dirac operator to both sides of our PDE to give that $\psi$ weakly solves
$$\Dl \psi = \Par^2 (\psi) = \Par (\Gamma_{\al} \gamma_{\al} \cdot \psi)$$ by Lichnerowitz's formula \cite{LM}.
Now, by extending $\psi$ and $\Gamma$ by zero and setting $V:= N[\Par (\Gamma_{\al} \gamma_{\al} \cdot \psi)]$ where $N$ is the Newtonian potential, we have that
$$|V|\leq I_1 [|\Gamma_{\al} \gamma_{\al} \cdot \psi|]$$where $I_1$ is the Riesz potential of order one - see \cite[p.3755-3756]{changyouwang_dirac}. By standard estimates (see e.g. \cite{gt}), we have
$$\|V\|_{L^{4}} \leq C\|\Gamma\|_{L^2}\|\psi\|_{L^{4}}.$$
Now, setting $k:=V-\psi$ we have that $\Dl k = 0$ in $B_1$ and $k\in L^4$ yielding (see \cite[Lemma 3.3.12]{helein_conservation})
$$\|k\|_{L^4(B_r)} \leq r^{\fr12} \|k\|_{L^{4}(B_1)} \leq  r^{\fr{1}{2}}(C\|\Gamma\|_{L^{2}}\|\psi\|_{L^{4}} + \|\psi\|_{L^{4}(B_1)}).$$
Therefore
\bee
\|\psi\|_{L^{4}(B_r)}^4 &\leq &(\|V\|_{L^{4}(B_1)} +\|k\|_{L^{4}(B_r)})^4\nn\\
&\leq & (C\|\Gamma\|_{L^{2}}\|\psi\|_{L^{4}(B_1)} + r^{\fr12}\|\psi\|_{L^{4}(B_1)})^4 \nn\\
&\leq & (1+\dl) r^2\|\psi\|_{L^{4}(B_1)}^4 + C_{\dl}\epsilon^4\|\psi\|_{L^{4}(B_1)}^4 \nn
\eee
by Young's inequality. Let $0<\nu <2$ and first set $\dl > 0$ and then $\epsilon >0$ small enough so that
$$\theta:=\fr{1+\dl + 4C_{\dl}\epsilon^4}{4} = \left(\fr12\right)^{\nu}.$$
Now by a re-scaling argument we have
$$\|\psi\|_{L^4(B_{\fr{R}{2}})}^4\leq \theta \|\psi\|_{L^4(B_{R})}^4.$$  Given $r<1$ we can find $k$ such that $$2^{-(k+1)}\leq r\leq 2^{-k}$$ and we have
$$\|\psi\|_{L^4(B_r)}^4 \leq \|\psi\|_{L^4(B_{2^{-k}})}^4 \leq \theta^k\|\psi\|_{L^4(B_1)}^4.$$
From here using standard techniques we can conclude that for all $2>\nu >0$ there exists $\epsilon > 0$ such that
$$\|\psi\|_{M^{4,2-\nu}(B_{\fr12})} \leq C\|\psi\|_{L^{4}(B_1)}.$$
Now this gives $\Gamma \psi \in M^{\fr43, \fr{2-\nu}{3}}_{loc}$. Setting $\nu:=\fr{16}{3q-4}$, and using the Riesz potential estimates of
Adams \cite{adams_riesz} gives $V\in M^{q, \fr{2q-8}{3q-4}}_{loc}$ yielding the appropriate estimate and eventually giving the proof.  \hfill$\square$

{\bf Proof of Theorem \ref{spinor_thm1}:}
First of all note that using \eqref{representation} and \eqref{diracoperator}, for $x\in B_1^+$ we have
\begin{eqnarray*}
\Par \psi (x)  &=& 2i\cvec{\pl{\psi_-}{z}(x)}{\pl{\psi_+}{\zb}(x)}\\
&=& \left[ \Gamma_1(x)\twomat0ii0 + \Gamma_2(x)\twomat0{1}{-1}0 \right]\cvec{\psi_+(x)}{\psi_-(x)}\\
&=& \cvec{(i\Gamma_1(x)+\Gamma_2(x))\psi_-(x)}{(i\Gamma_1(x) - \Gamma_2(x))\psi_+(x)}.
\end{eqnarray*}
 Extend the quantities arising in the theorem as follows:
$$\psi_E(x):=\twopartdef{\cvec{\psi_+(x)}{\psi_-(x)}}{x\in B_1^+}{\mp\cvec{R\psi_-(\rho(x))}{R\psi_+(\rho(x))}}{x\in B_1^-,}$$
$$\Gamma_E(x):=\twopartdef{\Gamma(x)}{x\in B_1^+}{R\rho^{\ast}\Gamma (x) R}{x\in B_1^-.}$$
We leave it to the reader to check that $\psi_E\in W^{1,1}(B_1)$ so we only have to check that $\psi_E$ solves the correct PDE on the lower half ball. Take $x\in B_1^-$ and check using \eqref{representation} and \eqref{diracoperator} that \begin{eqnarray*}
\Par \psi_E (x) &=& \mp 2i\cvec{R\pl{\psi_+}{\zb} (\rho(x))}{R\pl{\psi_-}{z} (\rho(x))} \\
&=& \mp\cvec{R[i\Gamma_1(\rho(x)) - \Gamma_2(\rho(x))]\psi_+(\rho(x))}{R[i\Gamma_1(\rho(x))  +\Gamma_2(\rho(x))] \psi_-(\rho(x))} \\
&=&\left[R\Gamma_1(\rho(x))R \twomat0ii0 -R\Gamma_2(\rho(x))R\twomat0{1}{-1}0\right] \cvec{\mp R\psi_- (\rho(x))}{\mp R\psi_+(\rho(x))} \\
&=& \Gamma_E(x)_{\al}\gamma_{\al} \cdot \psi_E (x)
\end{eqnarray*}
where we have used that $\Gamma = \Gamma_1 \ed x_1 + \Gamma_2 \ed x_2$ so that
$$R\rho^{\ast} \Gamma(x) R = R\Gamma_1(\rho(x))R \ed x_1 - R \Gamma_2(\rho(x)) R \ed x_2.$$
Therefore we have
\begin{equation*}
\Par \psi_E = (\Gamma_E)_{\al} \gamma_{\al}\cdot \psi_E\quad\text{in $B_1$}
\end{equation*}
at which point Theorem \ref{spinor_thm2} finishes the proof.              \hfill$\square$

\section{Regularity and smooth estimates}\label{highreg}
In this section, we first show the interior smooth estimates in all dimensions $m\geq2$,
Theorem \ref{thm: smooth estimates2d} and Theorem \ref{thm: smooth estimates}. Then, we prove the full regularity and smooth estimates at the free boundary
for weakly Dirac-harmonic maps in dimension $m=2$, Theorems \ref{thm: breg} and Theorem \ref{newthm}.

\subsection{Interior estimates in all dimensions, the proof of Theorems \ref{thm: smooth estimates2d} and \ref{thm: smooth estimates}}

We begin by proving the two dimensional result Theorem \ref{thm: smooth estimates2d} by first showing that the equation solved by the spinor, \eqref{spinor} is sub-critical. Then a simple application of  \cite[Theorem 1.1]{Sh_To} - or Theorem \ref{thm: main2} for $A-Id\equiv \zeta \equiv 0$ and $m=2$ - puts us in the position of being able to apply a simple bootstrapping argument.

{\bf Proof of Theorem \ref{thm: smooth estimates2d}:}
Looking back at \eqref{spinor} and setting $\Gamma^i_{j,\al} = \sum_{k,L} \pl{\nu_L^k}{z^j} \nu^i_L \phi^k_{\al}$ we have $|\Gamma (x)|\leq C|\D \phi (x)|$
 and therefore, we fix $q>4$ and conclude that if $\|\D \phi\|_{L^2}$ is sufficiently small then $\psi\in L^q_{loc}(B_1)$ with the appropriate estimate from Theorem \ref{spinor_thm2}. Next consider equation \eqref{map2} with $(g_1)_{ij} = \dl_{ij}$; we have $f\in L^p$ for some $1<p<2$ and
$$\|f\|_{L^p(B_{\fr34})}\leq C\|\D \phi\|_{L^2(B_1)}(1 + \|\psi\|_{L^4(B_1)}).$$
Now, applying \cite[Theorem 1.1]{Sh_To}, or indeed Theorem \ref{thm: main2} (for $m=2$, $\Om = \om$, $A\equiv Id$, $\zeta \equiv 0$) and
noting that $|\om(x)|\leq C|\D \phi(x)|$ tells us that when $\|\D \phi\|_{L^2(B_1)}$ is sufficiently small then we have (see remark \ref{rmkthmmain2})
$$\|\phi\|_{W^{2,p}(B_{\fr12})} \leq C\|\D \phi\|_{L^2(B_1)}(1 + \|\psi\|_{L^4(B_1)}).$$
We leave the rest of the proof to the reader since it is a simple bootstrapping argument using \eqref{map} and \eqref{spinor}.\hfill $\square$

{\bf Proof of Theorem \ref{thm: smooth estimates}:}
Consider now equation \eqref{map1} for $(g_1)_{ij} = \dl_{ij}$ and notice that we have $|\Om (x)|\leq C(|\D \phi(x)| + |\psi(x)|^2).$
Thus we can apply \cite[Corollary 1.3]{Sh}, also in \cite[Corollary 6.2.2]{thesis} or equivalently Theorem \ref{thm: main2} for
$A-Id\equiv \zeta \equiv f \equiv 0$ giving that (fixing $q>4$) there exists an $\epsilon>0$ such that if
$\|\D \phi\|_{M^{2,m-2}(B_1)} + \|\psi\|_{M^{4,m-2}(B_1)} \leq \epsilon$ then $\D \phi \in M^{q,m-2}_{loc}(B_1)$
with the appropriate estimate - see remark \ref{rmkthmmain2} for some details if necessary. We give a brief outline of the rest of the proof since it is essentially a standard boot-strapping argument but with some perhaps non-standard spaces.

Since $\D \phi \in M^{q,m-2}_{loc}(B_1)$, checking \eqref{spinor}, coupled with the fact that $\psi \in M^{4,m-2}$ gives $\Par \psi \in M_{loc}^{s,m-2}$ for $\fr1s = \fr14 + \fr1q$ and $s>2$ since $q>4$. Thus we have $\D \psi \in M^{s,m-2}_{loc}$ since Morrey spaces are preserved under singular integral transformations - see \cite{peetre_cam}. Now, one can easily use the Poincar\'e inequality and some standard facts about Morrey-Campanato spaces (see \cite{giaquinta} or \cite[p.33-34]{thesis}) to conclude $\psi\in C^{0,\fr{s-2}{s}}_{loc}$ with an estimate.
Looking back at \eqref{map} we see that the worst term on the right hand side is quadratic in $\D \phi$, which we can control in $M^{\fr{q}{2},m-2}_{loc}(B_1)$, therefore we have $\D^2 \phi \in M^{\fr{q}{2},m-2}_{loc}(B_1)$. Since $q>4$ and using a similar argument as above we can conclude $\D \phi \in C^{0,\fr{q-4}{q}}_{loc}(B_1)$. We now  leave the rest of the details to the reader since they follow from a straightforward bootstrapping argument using \eqref{map} and \eqref{spinor}.
\hfill$\square$

When $m=2$, since the problem is conformally invariant we can use Theorem \ref{thm: smooth estimates2d} to conclude the following:

\begin{cor}\label{cor_smooth_2d}
Let $(M^2, g)$  be a closed spin Riemann surface and $(N^d,h)$ a closed Riemannian manifold with $d\geq 2$.
Let $(\phi, \psi)$ be a weakly Dirac-harmonic map from $(M,g)$ to $(N,h)$, then it is smooth. Moreover for any $k\in \mathbb{N}$,
there exist $\epsilon = \epsilon (M,N)>0$ and $C=C(M,N,k)>0$ such that whenever
$$\int_{B^M_r(x_0)} |\D \phi|^2 \ed V_g  \leq \epsilon^2$$for some $x_0\in M$ and $0<r \leq \fr{i_M}{2}$, then
$$\|\D^k \phi\|_{L^{\infty}(B_{\fr{r}2}(x_0))} + r^{\fr12}\|\D^k \psi\|_{L^{\infty}(B_{\fr{r}2}(x_0))} \leq Cr^{-k}\|\D \phi\|_{L^{2}(B_r(x_0))}(1 + \|\psi\|_{L^{4}(B_r(x_0))}).$$
Here we have denoted intrinsic geodesic balls by $B^M$.
\end{cor}

We can couple Theorem \ref{thm: smooth estimates} with Remark \ref{remark general domain2} and recover the following result
(see \cite[Theorem 1.4]{WX}):
\begin{cor}
Let $(M^m, g)$ be a closed spin Riemannian manifold and $(N^d,h)$ a closed Riemannian manifold with $m,d\geq 2$. For any $k\in \mathbb{N}$,
there exist $\epsilon = \epsilon (M,N)>0$ and $C=C(M,N,k)>0$ such that if $(\phi, \psi)$ is a weakly Dirac-harmonic map from $(M,g)$ to  $(N,h)$ and
for some $x_0\in M$ and $0<r_0 \leq \fr{i_M}{2}$ there holds
$$\sup_{x\in B^M_{r_0}(x_0), \,\,0<r\leq r_0} \left\{ \fr{1}{r^{m-2}}\int_{B^M_r(x)} (|\D \phi|^2 + |\psi|^4)\ed V_g \right\} \leq \epsilon^2$$
then $(\phi, \psi)$ is smooth in $B^M_{r_0}(x_0)$ and the following estimates hold:
$$\|\D^k \phi\|_{L^{\infty}(B_{\fr{r}2}(x_0))} + r^{\fr12}\|\D^k \psi\|_{L^{\infty}(B_{\fr{r}2}(x_0))} \leq Cr^{-k}(\|\D \phi\|_{M^{2,m-2}(B_r(x_0))} + \|\psi\|_{M^{4,m-2}(B_r(x_0))}).$$\end{cor}

\subsection{Localisation at the free boundary in the domain and  target}\label{localisation}

Let $M$ be a compact spin Riemann surface and let $(\phi,\psi)\in \mathcal{X}^{1,2}_{1,\fr43}(M, N;\mathcal{S})$ be a weakly Dirac-harmonic map with free boundary $\phi(\partial M)$ on $\mathcal{S}$.
By conformal invariance in dimension $m=2$, for simplicity, we shall locate our problem in  a small neighbourhood of a boundary point and consider the case
that the domain is $B^+_1:= \left \{(x_1,x_2) \in \mathbb{R}^2 | x_1^2+x_2^2 < 1, x_2 \geq 0  \right \}$ and the free boundary
portion is $I:= \left \{(x_1,0) \in \mathbb{R}^2 | -1 < x_1 <  1 \right \}$. Moreover, we take $\g_{\a}= \partial x_{\a}, \a=1,2$.

It turns out that one can also localise the problem in the target. To see this, we give the following lemma, which is an improved version of Lemma 3.1 in \cite{CJWZ},
showing that the image of $\phi$ is contained in a small neighbourhood of some point $q \in \mathcal{S}$ under a smallness condition.

\begin{lemma} \label{nhdlemma} Let $N$ be a compact Riemannian manifold, isometrically
embedded in $\mathbb{R}^n$ and $\mathcal{S}\In N$ a closed submanifold. Then there is an $\epsilon_{0}=\epsilon_{0}(N)>0$ such that
if $(\phi,\psi)\in \mathcal{X}^{1,2}_{1,\fr43}(B_1^+, N;\mathcal{S})$ is a weakly Dirac-harmonic map with free boundary $\phi(I)$ on $\mathcal{S}$ satisfying
$$\|\D \phi\|_{L^2(B_1^+)} \leq \epsilon_0$$
then $${\rm dist}(\phi(x),\mathcal{S}) \leq C\epsilon_0(1+\|\psi\|_{L^4(B_1^+)}) \quad \text{for all $x\in B_{1/4}^+$}$$
with a constant $C=C(N)>0$. Moreover if we assume furthermore that $\|\psi\|_{L^4(B_1^+)}\leq \epsilon_0$, there
is a $q \in \mathcal{S}$ such that $\phi(x)\in B_{\ti{C}\epsilon_0}(q)$ for all $x\in B_{1/8}^+$ with $\ti{C}:=\ti{C}(N).$
\end{lemma}

{\bf Proof of Lemma \ref{nhdlemma}:}
The first statement follows by the proof of Lemma 3.1 in \cite{CJWZ} (see also \cite[Lemma 3.1]{Sc}). As such, it is sufficient to prove that given any $x_0 \in B_{\fr14}^+\sm I$ and $R:=\fr13 \text{dist}(x_0,I)$ then there exists a $C=C(N)>0$ with
$$|\D \phi (x)| \leq \fr{C(1+\|\psi\|_{L^4(B_1^+)})\epsilon_0}{R}$$
for all $x\in B_{2R}(x_0)$. This follows simply by applying the interior estimate from Corollary \ref{cor_smooth_2d} which gives that there exists a $C=C(N)>0$ with
$$|\D \phi (x)| \leq \fr{C\|\D \phi\|_{L^2(B_{3R}(x_0)}(1+\|\psi\|_{L^4(B_{3R}(x_0))})}{R}$$ for all such $x$.

We now prove the second statement:

Apply the $\epsilon$-regularity Theorem \ref{spinor_thm1} to the spinor equation \eqref{spinor} with boundary condition \eqref{B1}, we know that
$\|\psi\|_{W^{1,q}(B^+_{1/2})}\leq \|\psi\|_{L^4(B^+_1)}$ for any $1<q<2$ and hence by the boundary condition for the map \eqref{B2}, we have
$ (\pl{\phi}{\overrightarrow{n}})^{\top} =
{\rm Re}(  \mathcal{P}_{\mathcal{S}}(\overrightarrow{n}\cdot \psi^{\bot};\psi^{\top})) \in W^{1,p}_{\delta}(I_{1/2},\R^n)$ for any $1<p<2$,
where $I_{1/2}=I \ \cap B_{1/2}$. By Theorem \ref{Gn}, there is a $G \in W^{2,p}(B^+_{1/2}, \R^n)$ satisfying 
$$\|G\|_{W^{2,p}(B^+_{1/2})} \leq C (\|\nabla \phi\|_{L^{2}(B^+_1)}+\|\psi\|_{L^{4}(B^+_1)})$$ such that $ \pl{G}{\overrightarrow{n}}  =  (\pl{\phi}{\overrightarrow{n}})^{\top}$ on $I_{1/2}$ with $\overrightarrow{n}=-\gamma_2$, where $\top$ and  $\bot$ are with respect to $\mathcal{S}$. Then by definition of weakly Dirac-harmonic map with free boundary $\phi(I)$ on $\mathcal{S}$ and Green's identity, we have
\bee &&\int_{B^+_{1/4}} d \phi  \cdot \nabla V + \int_{B^+_{1/4}} \left \la \psi^i,
\g_{\a}\cdot \psi^j \right \ra  \left \la V_i, R(\phi)\left
(V,\phi_*(\g_\a)\right ) V_j \right \ra \nn\\
&=& \int_{I_{1/4}} (u_{(-\gamma_2)})^{\top} \cdot V, \nn\\
&=& \int_{B^+_{1/4}} d G \cdot \nabla^{\R^n} V  + \int_{B^+_{1/4}} \Delta G \cdot  V, \nn\\
&=& \int_{B_{1/4}^+} dG \cdot \nabla V + \int_{B_{1/4}^+} dG \cdot A(\nabla \phi, V)+ \int_{B_{1/4}^+}  \Delta G \cdot  V
\eee
for all compactly supported   $V\in H^1\cap L^{\infty}(B^+_{1/4}, \phi^{-1}TN)$ such that $V(x)
\in T_{\phi(x)}\mathcal{S}$ for $a.e. x \in I_{1/4}$.

As in \cite{CJWZ}, we use the geodesic reflection $\sigma$ to extend the two fields $(\phi,\psi)$ to the whole $B_{1/4}$. By Theorem \ref{Gn} again,
there is a $\widetilde{G} \in W^{2,p}(B^-_{1/4}, \R^n)$ such that $ d\widetilde{G}(\gamma_2)  =  ( D \sigma d \phi (\gamma_2))^{\top} =
(d \phi (-\gamma_2))^{\top}= {\rm Re}(  \mathcal{P}_{\mathcal{S}}(\overrightarrow{n}\cdot \psi^{\bot};\psi^{\top}))$ on $I_{1/4}$.
Apply similar arguments as in the proof of Theorem 3.1 in \cite{CJWZ} and apply Green's identity, we have
\bee \label{extended-eq}
&&\int_{B_{1/4}} d \phi \cdot_h \nabla^h V + \int_{B_{1/4}} \left \la
\psi^i, \g_{\a}\cdot \psi^j \right \ra  \left \la V_i,
R^h(\phi)\left (V,\phi_*(\partial x_\a)\right ) V_j \right \ra_h \nn \\
&=& \int_{I_{1/4}} (d\phi(-\gamma_2))^{\top} \cdot V + \int_{I_{1/4}} ((D \sigma d\phi)(\gamma_2))^{\top} \cdot D \sigma V   , \nn\\
&=& \int_{B_{1/4}^+} dG \cdot \nabla^{\R^n} V + \int_{B_{1/4}^+}  \Delta G \cdot  V + \int_{B_{1/4}^-} d\widetilde{G}  \cdot  \nabla^{\R^n} ( D\sigma V)
+ \int_{B_{1/4}^-} \Delta \widetilde{G} \cdot ( D\sigma V), \nn \\
&=& \int_{B_{1/4}^+} dG \cdot \nabla V + \int_{B_{1/4}^+} dG \cdot A(\nabla \phi, V)+ \int_{B_{1/4}^+}  (\Delta G) \cdot  V, \nn\\
&&+  \int_{B_{1/4}^-} d\widetilde{G}  \cdot  \nabla ( D\sigma V) +\int_{B_{1/4}^-} d\widetilde{G}  \cdot  A(\nabla \phi, D\sigma V)
+ \int_{B_{1/4}^-} (\Delta \widetilde{G}) \cdot ( D\sigma V), \nn \\
&=& \int_{B_{1/4}} H \cdot_{h} \nabla^{h} V + \int_{B_{1/4}} F \cdot_{h}  V,
\eee
for all compactly supported $V\in H^1\cap L^{\infty}(B_{1/4}, \phi^{-1}TN)$. Here
\bee \label{H}
H(x)=\left\{
\begin{array}{rcl}
(dG)^{\top},  \quad x \in B^+_{1/4}  \\
D \sigma (d \widetilde{G})^{\top},  \quad x \in B^-_{1/4}  \\
 \end{array}
\right.
\eee
and
\begin{equation} \label{F}
F(x)=\left\{
\begin{array}{rcl}
(\Delta G)^{\top} + (dG \cdot A(\nabla \phi, \cdot))^{\top},  \quad x \in B^+_{1/4}  \\
D \sigma (\Delta \widetilde{G})^{\top}+ D\sigma(d\widetilde{G}  \cdot  A(\nabla \phi,\cdot))^{\top},  \quad x \in B^-_{1/4}  \\
 \end{array}
\right.
\end{equation}
Here in \eqref{H} and \eqref{F}, $\top$ and $\bot$ are with respect to $N\subset \R^n$. It is easy to check that $H\in L^q$ for $q>2$ and $F \in L^{p}$ for $1<p<2$.

As in \cite{CJWZ}, let $e_i(x)\in T_{\phi} N$ be the $H^1$-tangent frame that is orthonormal with respect to the metric $h$,
for any $\varphi\in C^{\infty}_{0}(B_{1/4})$, taking $V=\varphi e_i$ in \eqref{extended-eq}, we derive the following equation:
\bee   d^*(d\phi \cdot_h e_i) = \left(  (\nabla^h e_i \cdot_h e_j)+ \mathfrak{R}_{ij} \right) (d\phi \cdot_h e_j) +  d^*(H_1) + F_1 . \eee
with $\mathfrak{R}_{ij}$ defined as in \cite{CJWZ} (section 3), $H_1\in L^q$ for $q>2$ and $F_1 \in L^{p}$ for $1<p<2$ . Then by slightly modifying the argument in the proof of Lemma 3.5 in \cite{CJWZ}, which is similar to \cite{riviere_struwe, WX}, we have (see Remark \ref{remark general domain})
$$ [\phi]_{C^{0,\a}(B_{1/8}^+)} \leq C (\parallel \nabla \phi \parallel_{L^2(B_1^+)}+||\psi||_{L^4(B_1^+)}).$$
The second statement follows immediately. \hfill$\square$

\subsection{Regularity and smooth estimates at the free boundary for weakly Dirac-harmonic maps from surfaces} \label{reg-smooth}

By conformal invariance in dimension $m=2$ and setting $\epsilon_0$ in Lemma \ref{nhdlemma} sufficiently small, we can without loss of generality assume that
$\phi(B^+_1) \subset \mathbf{U}$, where  $\mathbf{U}$ is a ball in $N$ centred at a point $q\in \mathcal{S}$  and
with Fermi coordinates $\{y^i\}_{i=1}^d$ about $\mathcal{S}$, therefore from now on we think of $\mathbf{U}\In \R^d$ as being some bounded domain.
Note that we are now in the set-up of Proposition \ref{pro3.7}.

We shall use an orthonormal frame to express the terms appearing there: Recall that the matrix $(g_{ij})=( \la \partial_i, \partial_j \ra_g )$ has block form
(see \eqref{Rg})
$ g =  \left(
           \begin{array}{cc}
             g_{ab} & 0\\
            0 & g_{\lambda \mu} \\
           \end{array}
         \right). $
Applying Gram-Schmidt orthonormalisation to the two frames  $\{\partial_j\}_{j=1}^s$ and $\{\partial_\mu\}_{\mu=s+1}^d$,
we obtain two smooth orthonormal frames  $\widehat{e}_{a}=\widehat{A}_{a}^{b} \partial_{b}$ and
$\widehat{e}_{\lambda}=\widehat{A}_{\lambda}^{\mu} \partial_{\mu}$, respectively.
Here $(\widehat{A}^{b}_{a}) \in C^{\infty}(\mathbf{U}, GL(s))$ and $(\widehat{A}_{\lambda}^{\mu}) \in C^{\infty}(\mathbf{U}, GL(d-s))$.
Clearly we have that $(\widehat{e}_{a}, \widehat{e}_{\lambda})_{1\leq a\leq s,  s+1\leq \lambda \leq d}$ form a  smooth orthonormal frame with respect to $g$.
Let now $$\widehat{\Om}^i_j:= ( \widehat{e}_i, \D^{N} \widehat{e}_j)$$ be the connection form for the Levi Civita connection on $TN$ over $\mathbf{U}$
with respect to the frame $\{\hat{e}_i\}$, noting that it is antisymmetric.

Denote $ \widehat{A} :=  \left(
           \begin{array}{cc}
            \widehat{A}^{a}_{b} & 0\\
            0 & \widehat{A}^{\lambda}_{\mu} \\
           \end{array}
         \right)  $
and set $A(x)=(\widehat{A}(\phi(x)))^{-1}$, then $A\in W^{1,2}\cap L^{\infty}(B_1^+, GL(d))$ and $A^{-1}\in W^{1,2}\cap L^{\infty}(B_1^+, GL(d))$.
It is easy to see that $A$ commutes with $R$ and there exist $\Lambda = \Lambda (N,\mathcal{S})>0$ and $C=C(N,\mathcal{S})>0$ such that
$$\Lambda^{-1}|\xi|\leq |A(x)\xi|\leq \Lambda |\xi|,   \quad |\D A(x)|\leq C|\D \phi(x)|$$ for almost every $x\in B_1^+$ and for all $\xi \in \R^d$.

Set $e_i(x)=\widehat{e}_i(\phi(x))$, $i=1,2,...,d$ and $\Om^i_j:= \phi^{\ast}\widehat{\Om}^i_j \in L^2(B_1^+,so(d)\otimes \wedge^1 \R^2)$
then $\{e_i\}_{i=1}^d\In W^{1,2}\cap L^{\infty}(B_1^+,\R^d)$ is an orthonormal frame with respective connection forms $\Om$
for $\phi^{-1}TN=B_1^+\times \R^d$, and there exists $C=C(N,\mathcal{S})>0$ such that
$$|\Om (x)|\leq C |\D \phi(x)|, \quad \text{for}\ a.e. \ x\in B_1^+$$  Consider $\eta\in C_c^{\infty}(B_1^+\cup I,\R^d)$ satisfying
$$\eta^i = 0 \quad \text{on $I$, when $s+1\leq i\leq d$}$$
and set $V:=\eta^i e_i$: notice that $V$ is tangent to $\mathcal{S}$ along $I$ and is as in Proposition \ref{pro3.7}.
We can now write
$$\D^{\phi} V = \ed \eta^i e_i + \eta^j\Om^i_j e_i.$$
Moreover we can write
$$\partial y^i(\phi(x)) = A^j_i(x) e_j(x)$$and therefore we set
$$f^n:=- A^n_h g^{hi}\mathcal{R}_{ijkl} \phi_{\ast}^j(\g_{\a})\la\psi^k,\g_\a \cdot \psi^l\ra, \quad \frak{g}^n:= A^n_c g^{cd}\Gamma^\lambda_{ad}g_{\lambda\mu}{\rm Re} (\la \overrightarrow{n}\cdot \psi^\mu,\psi^a\ra)$$
and noting that $\ed \phi = A^k_l\ed \phi^l e_k$, we have
$$\int_{M} \la A^k_l\ed \phi^l e_k, \ed \eta^i e_i + \eta^j\Om^i_j e_i\ra -\int_{M} \la f^n e_n, \eta^i e_i \ra
=  \int_{\partial M} \la \frak{g}^n e_n, \eta^i e_i \ra$$and since $\Om^T = - \Om$ and $e_i$ is an orthonormal frame we can write
\begin{equation*}\int_M ( A\ed \phi,\ed \eta) = \int_M  (\la \Om,A\ed \phi\ra + f, \eta) + \int_{\partial M} (\frak{G},\eta)\end{equation*}
which holds for any such $\eta$. Notice that $\frak{g}^n = 0$ for $s+1\leq n \leq d$ by definition and the properties of $A$ - and
we have once again denoted $\frak{G}:= (\frak{g}^1,\dots\frak{g}^s,0,\dots 0)$. Moreover by our choice of coordinates we have
$$\phi^i = 0 \quad \text{on $I$, when $s+1\leq i \leq d$}.$$
Clearly now we would like to apply Theorem \ref{thm: mainPDE}, however we have $f\in L^1$ and $\frak{G}\in W_{\partial}^{1,1}(I)$
so we first improve the regularity of our spinor.

Note that $\psi$ solves

\begin{eqnarray}
 &\Par \psi^i= - \Gamma^i_{jk}(\phi)\partial_\alpha\phi^j
\g_\alpha\cdot\psi^k,\quad i=1,2,\cdots,d,& \label{sp}\\
&\mathbf{B}_{\phi}^{\pm} \psi = 0\quad \text{on $I$}\label{sp1}.&
\end{eqnarray}
Therefore setting $(\Gamma^i_j)_\al = - \Gamma^i_{jk}(\phi)\phi^k_{\al}$ we can apply Theorem \ref{spinor_thm1} to conclude that for any $\fr43<r<2$ we have $\psi\in W^{1,r}_{loc}(B_1^+\cup I)$ and  there exists an $\epsilon=\epsilon(r)>0$ such that whenever $\|\D \phi\|_{L^2(B_1^+)} \leq \epsilon$ then
$$\|\D \psi\|_{L^r(B_{\fr78}^+)}\leq C\|\psi\|_{L^4(B_1^+)}.$$

By shrinking the domain slightly and for $1<p<2$, we simultaneously have $f\in L^p(B_1^+)$ and $\frak{G}\in W_{\partial}^{1,p}(I)$. Thus, applying Theorem \ref{thm: mainPDE} we can conclude that $\phi\in W^{2,p}_{loc}(B_1^+\cup I)$ with an appropriate estimate.

Now, the coupled system we are looking at is \eqref{sp} with boundary conditions \eqref{sp1}
and
\begin{eqnarray}
& -\Delta\phi ^{k}= \Gamma^k_{ij}(\phi)\phi_{\alpha}^{i}\phi_{\alpha}^{j}
  - \frac{1}{2}\mathcal{R}^k_{\ lij}(\phi) \la\psi^i,\nabla\phi^l\cdot\psi^j \ra, \quad
k=1,2,\cdots,d,&\nn\\
&\phi^k =0 \quad \text{on $I$,  for $s+1\leq k \leq d$}&\nn \\
&\pl{\phi^k}{\overrightarrow{n}} = A^k_c g^{cd}\Gamma^\lambda_{ad}g_{\lambda\mu}{\rm Re} (\la \overrightarrow{n}\cdot \psi^\mu,\psi^a\ra) \quad \text{on $I$, for $1\leq k\leq s$}\nn &
\end{eqnarray}
We now assume the following:
\begin{thm}\label{spinor_thm}
Let $k\in \mathbb{N}_0$, $1<p<\infty$, $U\In\R^2$ any open set with $T\In \partial U$ a smooth boundary portion. Suppose
$\psi =\cvec{\psi_+}{\psi_-} \in W^{1,p}(U,\mathbb{C}^2\otimes \R^d)$ solves the following boundary value problem:
\begin{equation*}
\Par \psi = F \in W^{k,p}(U, \mathbb{C}^2\otimes \R^d)\quad\text{in $U$}
\end{equation*} and
\begin{equation*}
\mathbf{B}_R^{\pm}\psi = 0 \quad \text{on $T\In\partial U$}\quad \iff \quad \psi_+ = \mp R\psi_- \quad \text{on $T\In\partial U$},
\end{equation*}
where $R$ is as in \eqref{R}. Then for any $V\cemb U\cup T$, $\psi \in W^{k+1,p}(V)$ and there exists some $C=C(p,V,T)>0$ such that
\begin{equation*}
\|\psi\|_{W^{k+1,p}(V)} \leq C(\|F\|_{W^{k,p}(U)}+\|\psi\|_{L^{p}(U)}).
\end{equation*}

\end{thm}

For us we have that $U=B_1^+$, $T=I$ and given the discussion in section \ref{recapfbvp} - see \eqref{localB1} - our spinor $\psi$ satisfies the conditions of the theorem.
A bootstrapping argument by repeatedly applying Theorem \ref{spinor_thm} with the results from appendix \ref{app: dl} finishes the proof of Theorem \ref{thm: breg}. \hfill $\square$

\begin{rmk}
We point out here that the chirality boundary condition for $\Par$ is (morally speaking) a traditional boundary condition for the Cauchy-Riemann operator $\db$. It is easy to check (and we will do so) that we can frame this system as a coupled system of $\db$ equations with a vanishing imaginary (or real) part on the boundary. Thus Theorem \ref{spinor_thm} can be thought of as classical.
\end{rmk}

\noindent{\bf Proof of Theorem \ref{spinor_thm}:}
We can de-couple the PDE, letting $F=\cvec{F_+}{F_-}$ we have that (see \eqref{diracoperator})
$$\Par \cvec{\psi_+}{\psi_-} = 2i \cvec{\partial \psi_-}{\overline\partial \psi_+}.$$
Therefore
$$\db (\psi_+) = -\fr{i}{2}F_-\quad\text{and}\quad  \db (\overline{\psi}_-) =  \fr{i}{2} \overline{F}_+, \quad \text{in $U$}$$
with $$\psi_+ \pm R\psi_- = 0\quad \text{on $\partial U$}.$$

Letting $h_1 : = \psi_+ \pm R\overline{\psi}_- \ \text{and} \
h_2 := \psi_+ \mp R\overline{\psi}_-$, we have that $h_j$ solve
$$\db (h_1) =-\fr{i}{2}( F_- \mp R \overline{F}_+)\quad \text{in $U$}, \quad Re(h_1) = 0 \quad \text{on $T\In\partial U$}$$ and
$$\db (h_2) = -\fr{i}{2}(F_- \pm R\overline{F}_+)\quad \text{in $U$}, \quad Im(h_2) = 0 \quad \text{on $T\In\partial U$}.$$
Now applying Theorem \ref{dbb} and noting that $\psi =\fr12\cvec{h_1 +h_2}{\pm R(\overline{h}_1 - \overline{h}_2)}$ we have completed the proof of Theorem \ref{spinor_thm}. \hfill $\square$




\appendix

\section{Morrey Spaces}\label{def and not}

We introduce the Morrey spaces $M^{p,\beta}(E)$ for $1\leq p<\infty$ and $0\leq\beta\leq m$ ($E \In \R^m$). We say that $g \in M^{p,\beta}(E)$ if $M_{\beta}[g^p](x):= \sup_{r>0}  r^{-\beta} \int_{B_r(x)\cap E} |g|^p \in L^{\infty}$ with norm (which makes $M^{p,\beta}$ a Banach space) $\|g\|_{M^{p,\beta}(E)}^p = \|M_{\beta}[g^p]\|_{L^{\infty}(E)}.$

We note the following extension theorem for Sobolev-Morrey spaces:

\begin{thm}
The unit ball $B_1\In \R^m$ is an extension domain for $M_1^{p,\beta}$. Given $v\in M_1^{p,\beta}(B_1)$ there exists $\ti{v}\in M_1^{p,\beta}(\R^m)$ such that
$$\|\ti{v}\|_{M^{p,\beta}} + \|\D \ti{v}\|_{M^{p,\beta}} \leq C(\|v\|_{M^{p,\beta}(B_1)} + \|\D v\|_{M^{p,\beta}(B_1)})$$
and $v=\tilde{v}$ almost everywhere.
\end{thm}
We could not find a proof of this theorem, however it follows easily by standard techniques for extension theorems. In general it should remain true for any domain $U$ with $\partial U$ as above but $\phi$ must be Lipschitz. 

We also note that if $\int_{B_1} v = 0$ then $$\|\D \ti{v}\|_{M^{p,\beta}} \leq C([v]_{\Ca^{p,\beta}(B_1)} + \|\D v\|_{M^{p,\beta}(B_1)}) \leq C\|\D v\|_{M^{p,\beta}(B_1)}$$
where $\Ca^{p,\beta}$ is the Campanato space - see \cite{giaquinta}.

\section{Classical boundary value estimates} \label{app: dl}
Here we recall results about the classical boundary value problems for the Laplacian and the Cauchy-Riemann operator, we refer the reader to
\cite{wehrheim, agmon_doug_niren} for background material. For all of the theorems below $U\In \R^m$ is any open domain and $T\In \partial U$ is a smooth boundary portion.

\begin{thm}\label{dlb}
Let $k\in \mathbb{N}_0$ and $1<p<\infty$. Suppose that  $u\in W^{1,p}$ weakly solves
\bee
-\Dl u  &=& f\in W^{k,p}(U), \\
u&=&0 \quad \text{on $T$}
\eee
then for any $V\cemb U\cup T$, $u\in W^{k+2, p}(V)$ and there exists some $C=C(p,k,V,T)>0$ such that
\begin{equation*}
 \|u\|_{W^{k+2,p}(V)} \leq C(\|f\|_{W^{k,p}(U)} +\|u\|_{L^p(U)}).\end{equation*}
\end{thm}

\begin{thm}\label{dlnb}
Let $k\in \mathbb{N}_0$ and $1<p<\infty$. Suppose that  $u\in W^{1,p}$ weakly solves
\bee
-\Dl u  &=& f\in W^{k,p}(U), \\
\pl{u}{\overrightarrow{n}} &=& g \in  W_{\partial}^{k+1,p}(T)
\eee
then for any $V \cemb U\cup T$, $u\in W^{k+2, p}(V)$ and there exists some $C=C(p,k,V,T)>0$ such that
\begin{equation*}
 \|u\|_{W^{k+2,p}(V)} \leq C(\|f\|_{W^{k,p}(U)}+\|g\|_{W_{\partial}^{k+1,p}(T)} +
 \|u\|_{L^p(U)}).\end{equation*}
\end{thm}

We also recall the analogue for the Cauchy-Riemann operator in $\mathbb{C}$.

\begin{thm}\label{dbb}
Let $U\In \mathbb{C}$ be any domain and $T\In \partial U$ a smooth boundary portion, $k\in \mathbb{N}_0$ and $1<p<\infty$. Suppose that $h\in W^{1,p}$ solves
\bee
\db h  = f\in W^{k,p}(U), \eee
\bee
Re(h) =0 \quad \text{or}\quad Im(h) = 0 \quad \text{on $T$}
\eee
then for any $V\cemb  U\cup T$, $h\in W^{k+1, p}(V)$ and there exists some $C=C(p,k,V,T)>0$ such that
\begin{equation*}
 \|h\|_{W^{k+1,p}(V)} \leq C(\|f\|_{W^{k,p}(U)} +\|h\|_{L^{p}(U)}).
\end{equation*}
\end{thm}

\end{document}